\documentclass[12pt]{article}
\usepackage{amssymb}
\usepackage{amsmath}
\usepackage{eucal}
\usepackage{graphicx}

\def \ind {\hbox{ 1\hskip -3pt I}}
\def \Z  {\mathbb{Z}}  
\def \N  {\mathbb{N}}
\def \E  {\mathbb{E}}
\def \P  {\mathbb{P}} 
\def \R  {\mathbb{R}} 
\def \G {{\mathcal G}}

\newcommand {\cro}[1] {\left[ {#1} \right]}
\newcommand {\acc}[1] {\left\{ {#1} \right\}}
\newcommand {\pare}[1] {\left( {#1} \right)}

\newtheorem{theo}{Theorem}[section]
\newtheorem{prop}[theo]{Proposition}
\newtheorem{lemma}[theo]{Lemma}

\newtheorem{rem}[theo]{Remark}
\newcommand {\refeq}[1] {(\ref{#1})}

\def \DD {{\mathcal D}}
\def \Db{\underline{{\mathcal D}}_b}
\def \Dh{{\bar{\mathcal D}}_b}
\def \RR {{\mathcal R}}

\def\reff#1{(\ref{#1})}
\begin{document}
%\Large
\title{Self-Intersection Times for Random Walk, \\
and Random Walk in Random Scenery\\
in dimensions $d\geq 5$.}
\author{Amine Asselah \& Fabienne Castell\\C.M.I., Universit\'e de Provence,\\
39 Rue Joliot-Curie, \\F-13453 Marseille cedex 13, France\\
asselah@cmi.univ-mrs.fr \& castell@cmi.univ-mrs.fr}
\date{}
\maketitle
\begin{abstract}
Let $\{S_k,k\geq 0\}$ be a symmetric random walk
on $\Z^d$, and $\{\eta(x),x\in \Z^d\}$ an independent random field
of centered i.i.d.\  with tail decay $P(\eta(x)>t)\approx\exp(-t^{\alpha})$.
We consider a Random Walk in Random Scenery, that is 
$X_n=\eta(S_0)+\dots+\eta(S_n)$.
We present asymptotics for the probability, over both randomness,
that $\{X_n>n^{\beta}\}$ for $\beta>1/2$ and $\alpha>1$.
To obtain such asymptotics, we establish large deviations estimates 
for the self-intersection local times process $\sum l_n^2(x)$, where $l_n(x)$
is the number of visits of site $x$ up to time $n$.
\end{abstract}



{\em Running}: Random Walk in Random Scenery.

%___________________________________________________________________
% SECTION 1: INTRODUCTION.
%____________________________________________________________________
\section{Introduction.}
\label{intro}
We study transport in divergence free
random velocity fields. For simplicity,
we discretize both space and time and consider the simplest
model of shear flow velocity fields:
\[
\forall x,y\in \Z\times \Z^d,\qquad V(x,y)=\eta(y) e_x,
\]
where $e_x$ is a unit vector in the first coordinate of $\Z^{d+1}$, and
$\{\eta(y),y\in \Z^d\}$ are i.i.d.\  real random variables.
Thus, space consists of the sites of the cubic lattice $\Z^{d+1}$ and
the direction of the shear flow is $e_x$. We wish to model a polluant evolving
by two mechanisms: a passive transport by the velocity field, and
collisions with the other fluid particles modeled by
random centered and independent
increments $\{(\alpha_n,\beta_n)\in \Z\times \Z^d,n\in \N\}$,
independent of the velocity field.
Thus, if $R_n\in \Z\times\Z^d$ is the polluant's position at time $n$, then
\begin{equation} 
\label{eq-intro.1}
R_{n+1}-R_n=V(R_n)+(\alpha_{n+1},\beta_{n+1}),\quad\text{and}\quad
R_0=(0,0).
\end{equation} 
When solving by induction for $R_n$, \reff{eq-intro.1} yields
\begin{equation}\label{eq-intro.2}
R_n=\pare{\sum_{k=1}^n \alpha_k+\eta(0)+\sum_{k=1}^n \eta(\sum_{i=1}^k \beta_i),
\sum_{k=1}^n \beta_k}.
\end{equation} 
The sum $\beta_1+\dots+\beta_n$ is denoted by $S_n$, and called the Random Walk (RW).
The displacement along $e_x$ consists of two independent parts: a sum
of i.i.d.\  random variables $\alpha_1+\dots+\alpha_n$, and
a sum of {\it dependent} random variables $\eta(S_0)+\dots+\eta(S_n)$, which
we denote by $X_n$ and call the Random Walk in Random Scenery (RWRS). 
Writing  it in terms of local times of the RW, say $\{l_n(x),\ n\in\N,\ x\in \Z^d\}$,
we get 
\begin{equation} 
\label{eq-intro.3}
X_n=\sum_{k=0}^n \eta(S_k)=\sum_{x\in \Z^d} l_n(x) \eta(x),\quad
\text{where} \quad l_n(x)=\sum_{k=0}^n \ind\{S_k=x\}.
\end{equation} 
The process $\{X_n,n\in \N\}$ was studied at about the same time
by Kesten \& Spitzer \cite{ks},  
Borodin \cite{B1, B2}, and Matheron \& de Marsily \cite{matheron-demarsily}.
The fact that in dimension 1,  $E[X^2_n]\sim n^{3/2}$ made the model
popular and led the way to examples of {\it superdiffusive} behaviour. 
However, the {\it typical} behaviour of $X_n$ resembles 
that of a sum of $n$ independent variables all the more when dimension is large.

Our goal is to {\it estimate} the {\it probability} that $X_n$ be {\it large}.
By {\it probability}, we consider averages with respect to the two randomness,
and $P=\P_0\otimes P_{\eta}$, where $\P_0$ is the law of
the nearest neighbors symmetric random walk $\{S_k,k\in \N\}$ on $\Z^d$ with $S_0=0$,
and $P_{\eta}$ is the law of the velocity field. 

Now, when $d\ge 3$, Kesten \& Spitzer established in~\cite{ks} 
that $X_n/\sqrt{n}$ converges in law to a Gaussian variable. Thus,
by {\it large}, we mean $\{X_n>n^{\beta}\}$ with $\beta>1/2$. We expect
$P(X_n>n^{\beta})\approx \exp(-In^\zeta )$ with constant rate $I>0$, 
and we characterize in this work the exponent $\zeta$.
For this purpose, the only important feature of the $\eta$-variables is
the $\alpha$-exponent in the tail decay:
\begin{equation}\label{eq-intro.4}
\lim_{t\to\infty} \frac{\log P_{\eta}(\eta(0)>t)}{t^{\alpha}}=-c,
\quad\text{for a positive constant  }c.
\end{equation} 

Let us now recall the classical estimates for $P(Y_1+\dots+Y_n>n^\beta)$, where
$\beta>1/2$ and the $\{Y_n,n\in \N\}$ are centered i.i.d.\  
with tail decay $P(Y_n>t)\approx\exp(-t^{a})$, with $a>0$.
There is a dichotomy between a
``collective'' and an ``extreme'' behaviour. In the former case,
each variable contributes about the same, whereas in latter case, only
one term exceeds the level $n^\beta$, when the others remain small.
Thus, it is well known that
$P(Y_1+\dots+Y_n>n^\beta)\sim \exp(- In^\zeta)$ with three regimes for
the exponent $\zeta$.
\begin{itemize}
\item When $\beta\ge 1$ and $a>1$, a large collective contribution yields
$\zeta=(\beta-1)a+1$.
\item When $\beta<1$ and $\beta(2-a)<1$, 
a small collective contribution yields $\zeta=2\beta-1$.
\item When  $\beta>1/(2-a)$ and $a<1$, an extreme contribution yields $\zeta=\beta a$.
\end{itemize}
\begin{figure}
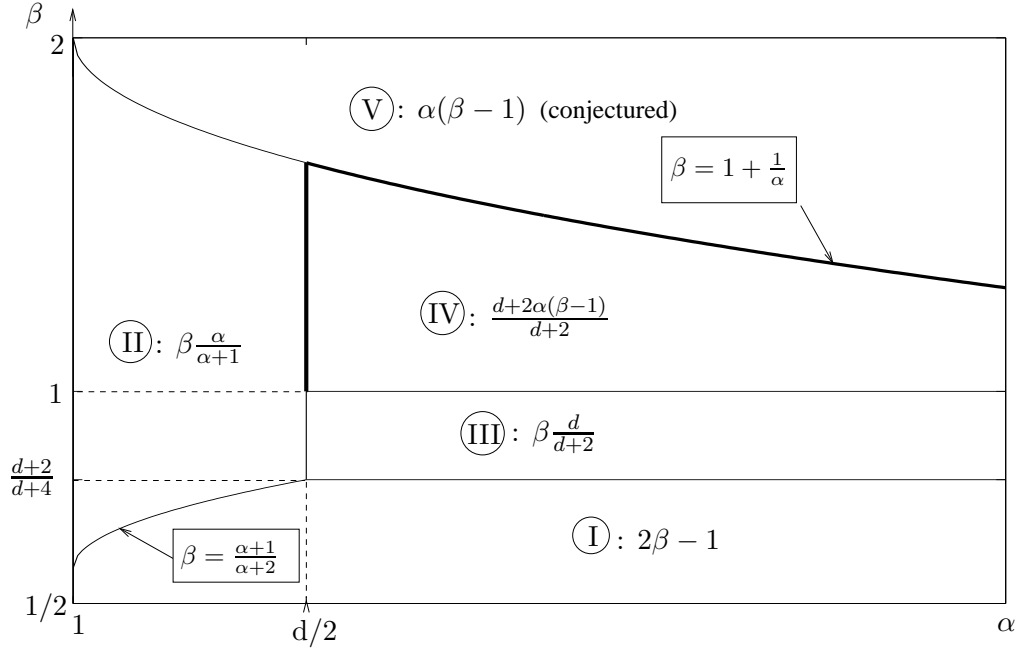

\input phase.pstex_t
\caption{$\zeta$-exponent diagram}
\label{diag-phase}
\end{figure}
For the RWRS, one expects a rich interplay between the scenery
and the random walk.
To get some intuition about the expression of $\zeta$ in terms
of $\alpha$ and $\beta$,
we propose simple scenarii leading to Figure~\ref{diag-phase}.
Here also, we focus on the exponent, and constants are omitted.

\begin{itemize} 
\item {\bf Region I}. No constraint is put on the walk. When $d\ge 3$, the range 
of the walk is of order $n$ and visited sites are typically visited once. Thus,
$\{X_n>n^\beta\} \sim \{\eta_1+\dots+\eta_n>n^\beta\}$. When $\beta<1$, the latter 
sum performs a moderate deviations of order $n^{\beta}$. Since the $\eta$-variables
satisfy Cramer's condition, we obtain $P(X_n>n^\beta) \ge\exp(-n^{2 \beta-1})$.
Thus, the $\zeta$-exponent in Region I is $\zeta_I=2\beta-1$.
\item {\bf Region II, V}. A few sites are visited often, 
so that $X_n\sim \eta(0)l_n(0)$.
Now, using the tail behaviour of $\eta(x)$, and the fact that in $d \geq 3$,
$l_n(x)$ is bounded by an exponential variable, we obtain
\begin{eqnarray*} 
 P\pare{X_n \geq n^{\beta}} 
& \ge& P\pare{l_n(0) \eta(0) \geq n^{\beta}}
\sim\sup_{k\le n}\pare{
 \P_0 \pare{l_n(0) =k} P_{\eta}\pare{\eta(0) \geq \frac{n^\beta}{k}}}\cr
& \sim & \exp\pare{-\inf_{k\le n}\pare{ k+(\frac{n^\beta}{k})^{\alpha}}}  \, . 
\end{eqnarray*}
Now, the minimum of $k\mapsto k+n^{\beta\alpha}/k^\alpha$ is reached for
$k^*=n^{\beta\alpha /(\alpha+1)}$. Since, we impose also that $k\le n$, two
different exponents prevail according to the value of $\beta$:
\begin{itemize}
\item[(II)] $\beta<(\alpha+1)/\alpha$, and $\zeta_{I\!\!I}:=\beta\alpha /(\alpha+1)
<1$. The RW spends a time of order $n^{\zeta_{I\!\!I}}$ on favorite sites.
\item[(V)] $\beta\ge (\alpha+1)/\alpha$, and $\zeta_{V}:=\alpha(\beta-1)$.
The RW spends a time of the order of $n$ on favorite sites.
\end{itemize}
\item {\bf Region III, IV}. The random walk is localized a time $T$
in a ball $B_r$
of radius $r$, with $r^2 \ll T$: this costs of the order of $\exp(-T/r^2)$. 
Then, during this period, each site of $B_r$ is visited about $T/r^d$, 
and we further assume that $r^d\ll T$. Thus
\begin{equation}\label{region3}
P\pare{X_n \geq n^{\beta}} \geq \exp(-\frac{T}{r^2}) 
P_{\eta} \pare{\frac{1}{\sqrt{r^d}} \sum_{B_r} \eta_j \geq 
\frac{n^{\beta} r^{d/2}}{T}} \, .
\end{equation} 
Two different exponents prevail according to $\beta$:
\begin{itemize}
\item[(III)] $\beta\le 1$. The condition $1\ll n^{\beta}r^{d/2}/T\le r^{d/2}$
means that the sum of $\eta$-s
performs a moderate (up to large) deviations and this costs of the order of
$\exp(-n^{2\beta}r^{d}/T^2)$. When the two costs are equalized and the
paramter $r$ and $T$ optimized,
we obtain that the walk is localized a time $T=n^\beta$ on a ball of volume
$r^d=n^{\zeta_{I\!\!I\!\!I}}$, with $\zeta_{I\!\!I\!\!I}:=d\beta/(d+2)$.
\item[(IV)] $\beta>1$. Here $T=n$ and we deal with a 
very large deviations for a sum of i.i.d.\ . This has a cost of order 
$\exp(-n^{\alpha(\beta -1)}r^d)$.
Choosing $r$ so that $n/r^2 = n^{\alpha(\beta -1)}r^d$, we obtain
$\zeta_{I\!V}:=(d+2\alpha(\beta-1))/(d+2)$. The condition $r \gg 1$ is equivalent 
to $\beta < 1 + 1/\alpha$. The walk is localized all the time
on a ball of radius $r$ satisfying $r^{d+2}=n^{1-\alpha(\beta-1)}$.
\end{itemize}
\end{itemize} 
 
The following regions have already been studied.
\begin{itemize}
\item $\alpha = +\infty$ (bounded scenery) and $\beta =1$ in \cite{AC03}
(actually Brownian motion is considered there instead of RW).
\item $\alpha =2$ (Gaussian scenery) and $\beta \in [1, 1+ 1/\alpha]$ 
in \cite{FC, CP}.
\item Region IV ($\alpha > d/2$, $1 \leq \beta < 1+ 1/\alpha$) in \cite{GKS}.
\item $0<\alpha<1$ and $\beta>\frac{1+\alpha}{2}$, in $d\ge 3$, in \cite{HGK}.
This region is outside Figure 1.
\item $\beta=1$ and $\alpha < d/2$ in \cite{AC04}. Contrary to the previous
cases, distinct lower and upper bounds with the same exponent are obtained 
in \cite{AC04}.
\end{itemize} 

This paper is devoted to regions I, II and III. Henceforth, we
consider $d\ge 5$, unless explicitly mentioned.

\begin{prop}\label{propMDUB} Upper Bounds for the RWRS.
\begin{enumerate} 
\item Region I. 
We assume $\beta < \min(\frac{\alpha + 1}{\alpha+2},
\frac{d/2+ 1}{d/2+2})$. There exists an explicit $y_0$, such that for
$y >y_0$, there exists a constant ${\bar c_1} > 0$ such that 
\begin{equation} \label{MDUBR1}
P(X_n\geq n^{\beta} y ) \leq \exp(- {\bar c_1}  n^{2\beta-1}).
\end{equation} 
\item Region II. Let $\alpha < d/2$, and $\beta \ge \frac{\alpha +1}{\alpha +2}$.
For $y > 0$, there exists a constant ${\bar c_2} > 0$, such that  
\begin{equation} 
\label{MDUBR2}
P(X_n\geq n^{\beta} y ) \leq  \exp(-{\bar c_2}  n^{\beta\alpha/(\alpha +1)}).
\end{equation} 
For the case $\beta=\frac{\alpha +1}{\alpha +2}$, we further assume that $y>y_0$.

Moreover, when $\beta >\frac{\alpha +1}{\alpha +2}$, the main contribution
to $\{X_n\geq n^{\beta} y \}$ comes from the level sets
\[
\DD_{I\!\!I}:=\{x:\ n^{b-\delta}<l_n(x)<n^{b+\delta}\}\quad\text{ with }
\quad b=\frac{\beta}{\alpha+1},\quad\text{and any }\delta>0.
\]
In other words, for any $y>0$
\begin{equation} \label{major-RegII}
\lim_{n\to\infty} \frac{1}{n^{\zeta_{I\!\!I}}}
\log P\pare{\sum_{x\not\in \DD_{I\!\!I}}\eta(x)l_n(x)>y n^\beta}=-\infty.
\end{equation}
\item Region III. Let $\alpha \geq d/2$ and $\beta \ge \frac{d/2+1}{d/2+2}$.
For $y > 0$ and $\epsilon >0$ small, there exists a constant ${\bar c_3} > 0$, 
such that  
\begin{equation} 
\label{MDUBR3}
P(X_n\geq n^{\beta} y ) \leq  \exp(-\bar{c_3}  n^{\beta \frac{d}{d+2} -\epsilon}).
\end{equation}   
For the case $\beta=\frac{d/2+1}{d/2+2}$, we further assume that $y>y_0$.

Moreover, if we define $\DD_{I\!\!I\!\!I}:=\{x:\ 0<l_n(x)<n^{b+\delta}\}$ with
$b:=\beta/(d/2+1)$, then we have , for any $\delta>0$ small enough
\begin{equation} \label{major-RegIII}
\lim_{n\to\infty} \frac{1}{n^{\zeta_{I\!\!I\!\!I}}}
\log P\pare{\sum_{x\not\in \DD_{I\!\!I\!\!I}}\eta(x)l_n(x)>y n^\beta}=-\infty.
\end{equation}

\end{enumerate}  
\end{prop} 
\begin{rem} Note that the control in Region III is less satisfactory than in Regions
I and II. An inspection of the proof makes it clear that our techniques actually
yield a logarithmic artifact $P(X_n\geq n^{\beta} y ) 
\leq \exp(-n^{\zeta_{I\!\!I\!\!I}}/
\log(n))$.
\end{rem}

We indicate below
lower bounds for $P(X_n \geq n^{\beta}y)$, which prove that we 
obtain the correct rates of the logarithmic decay of $P(X_n \geq n^{\beta}y)$.
These lower bounds are given under an additional symmetry assumption on
the scenery, which is not crucial, but simplifies the proofs. Hence,
we say that a real random variable is bell-shaped,  if its law has
a density with respect to Lebesgue which is even, and decreasing on $\R^+$.
\begin{prop}\label{propMDLB} Lower Bounds for the RWRS. \\
Assume $d\geq 3$, and that the random variables $\{\eta(x), x\in \Z^d\}$ 
are bell-shaped.
\begin{enumerate} 
\item Region I. Let $1 \geq \beta > 1/2$. For any $y > 0$, there exists 
a constant $\underline{c}_1 > 0$, such that 
\begin{equation}\label{MDLBR1}
P(X_n\geq n^{\beta} y) \geq  \exp(- \underline{c}_1 n^{2 \beta-1}) \, .
\end{equation} 
\item  Region II. Let  $\beta \leq 1 + 1/\alpha$. For any $y > 0$, there exists 
a constant $\underline{c}_2 > 0$, such that 
\begin{equation} 
\label{MDLBR2}
P(X_n\geq n^{\beta} y)\geq  \exp(- \underline{c}_2 n^{\beta\alpha/(\alpha +1)}).
\end{equation} 
\item  Region III. Let  $\beta \le 1$. For any $y > 0$, there exists 
a constant $\underline{c}_3 > 0$, such that 
\begin{equation} 
\label{MDLBR3}
P(X_n\geq n^{\beta} y)\geq  \exp(- \underline{c}_3 
n^{\beta\frac{d}{d+2}}).
\end{equation} 
\end{enumerate} 
\end{prop}

In the process of establishing Proposition\reff{propMDUB},
one faces the problem of evaluating the chances the random walk visits often
the same sites. More precisely, a crucial quantity is the 
{\it self-intersection local time process} (SILT):
\begin{equation} 
\label{eq-intro.7}
\Sigma_n^2=\sum_{x\in \Z^d} l_n^2(x) = n+1+ 2 \sum_{0\leq k<k'\leq n} 
\ind\{S_k = S_{k'}\}.
\end{equation} 

It is expected that $\Sigma^2_n$ would show up in the study of RWRS.
Indeed, $\Sigma_n^2$ is the variance of $X_n$ when averaged over $P_{\eta}$.
If we assume for a moment that the $\eta$-variables are standard Gaussian,
then conditionally on the random walk, $X_n$ is a Gaussian variable with variance
$\Sigma_n^2$, so that 
\begin{equation} 
\label {eq-intro.8}
P_{\eta}(\sum_{x\in \Z^d} \eta(x)l_n(x)>n^{\beta})\leq
\exp\left(-\frac{n^{2\beta}}{2\sum_{x\in \Z^d}l^2_n(x)}\right)
\end{equation} 
It is well known that an inequality similar to
\reff{eq-intro.8} holds for any tail behaviour
\reff{eq-intro.4} with $\alpha \geq 2$. Now, if we average with respect 
to the random walk law, then for any $\gamma>0$
\begin{eqnarray} 
\label{eq-intro.9}
P(X_n>n^{\beta})&\leq& 
\E_0\cro{\exp\left(-\frac{n^{2\beta}}{2\sum_{x\in \Z^d}l^2_n(x)}\right)}
\cr
&\leq & \exp(- n^{2\beta-\gamma})+
\P_0\left(\sum_{x\in \Z^d}l^2_n(x)> n^{\gamma}\right).
\end{eqnarray} 
Hence, at least for large $\alpha$, we have to evaluate the logarithmic
decay of quantities such as 
$\P_0\left(\sum_{x\in \Z^d}l^2_n(x)> n^{\gamma}\right)$. 
Note first that for $d \geq 3$, and $ n \rightarrow \infty$, 
\begin{equation} 
\label{moy-SILT.eq}
\E_0\cro{\sum_{x \in \Z^d}  l_n^2(x)} \simeq n (2 G_d(0)-1) \, ,
\end{equation}   
where  $G_d$ is the Green kernel 
\[
G_d(x) \triangleq \E_0\cro{l_{\infty}(x)} \, .
\]
Therefore, we have to take $\gamma \geq 1$ to be in a large deviations
scaling. 
For large deviations of SILT in $d=1$, we refer the reader
to Mansmann \cite{Ma91}, and Chen \& Li \cite{chen-li}, while in 
$d=2$, this problem is treated in Bass \& Chen \cite{bass-chen}, and
in Bass, Chen \& Rosen \cite{bass-chen-rosen}.
 
We first present large deviations estimates for the SILT. 
\begin{prop}\label{RLDSILT.prop} Assume $d\ge 5$.
For $y>1 + 2 \sum_{x \in \Z^d} G_d(x)^2$, there are positive
constants $\underline{c},\bar{c}$ such that
\begin{equation} 
\label{LDLBSILT.eq}
\exp(-\bar{c} \sqrt{n})\ge
\P_0 \pare{\, \sum_{x \in \Z^d} l_n^2(x)  \geq ny} 
\geq \exp(-\underline{c} \sqrt{n}) \, .
\end{equation}   
\end{prop}

Proposition~\ref{RLDSILT.prop} is a corollary of the next result
where we prove that the main contribution in the estimates 
comes from the region where the local time is of order $\sqrt{n}$. 

\begin{prop}\label{LDSILT.prop}
\begin{enumerate} 
\item For $\epsilon > 0$, and $y > 1 + 2 \sum_{x \in \Z^d} G_d(x)^2$, 
\begin{equation} 
\label{LDUB1.eq}
\limsup_{n \rightarrow \infty} \frac{1}{\sqrt{n}} \log 
\P_0 \pare{\sum_{x: l_n(x) \leq n^{1/2 - \epsilon}} l_n^2(x) \geq ny} 
= - \infty \, .
\end{equation} 
\item For $y>0$ and $\epsilon>0$, there exists a constant $\tilde{c} > 0$, such that 
\begin{equation} 
\label{LDUB2.eq} 
\limsup_{n \rightarrow \infty} \frac{1}{\sqrt{n}} \log 
\P_0 \pare{\sum_{x: l_n(x) >  n^{1/2 - \epsilon}} l_n^2(x) \geq ny}
\leq - \tilde{c}  \, .
\end{equation} 
\end{enumerate} 
\end{prop} 
We present now estimates for $\P_0(\sum l^p_n(x)>n^\gamma)$.
\begin{prop}\label{lnp-prop}
(i) Assume $p>\frac{d}{d-2}$, and $p>\gamma>1+\frac{(p-2)^+}{(d-2)p+4}$.
There are $c_1,c_2>0$ such that
\begin{equation}
\label{lnp-eq1}
e^{-c_1n^{\gamma/p}}\le \P_0\pare{\sum_{x\in \Z^d} l^p_n(x)>n^\gamma}
\le e^{-c_2n^{\gamma/p}}.
\end{equation}
(ii) Assume $1<p\le\frac{d}{d-2}$, and $p>\gamma>1$. For any $\epsilon>0$,
there are $d_1,d_2>0$ such that
\begin{equation}
\label{lnp-eq2}
e^{-d_1n^{\zeta}}\le \P_0\pare{\sum_{x\in \Z^d}\!\! l^p_n(x)>n^\gamma}
\le e^{-d_2n^{\zeta-\epsilon}},\quad\text{with}\quad 
\zeta=1-\frac{2}{d}\frac{(p-\gamma)}{(p-1)}.
\end{equation}
\end{prop}
Let us give some heuristics on the proof of 
Proposition \ref{LDSILT.prop}. First of all, we decompose $\Sigma_n^2$ using the 
level sets of the local time. Note that it is not useful to consider $\{x:l_n(x)\gg
\sqrt{n}\}$, since $l_n(x)$ is bounded by an exponential variable. Now,
for a subdivision $\{b_i\}_{i \in \N}$ of $[0,1/2]$,
let $\DD_{b_i} = \acc{x \in \Z^d:\ n^{b_i} \leq l_n(x) < n^{b_{i+1}}}$.
Denoting by $|\Lambda|$ the number of sites in $\Lambda \subset \Z^d$, we then have 
\[\Sigma_n^2 = \sum_i \sum_{x \in \DD_{b_i}} l^2_n(x) 
\leq \sum_i n^{2b_{i+1}} |\DD_{b_i}|\, .
\] 
Hence, choosing $(y_{b_i})_{i \in \N}$ such that $\sum_i y_{b_i} \leq y$,
\[
\P_0\pare{\Sigma_n^2 \geq n y} 
\leq \sum_i \P_0 \pare{\sum_{x \in \DD_{b_i}} l_n^2(x) \geq n y_{b_i}}
\leq \sum_i \P_0 \pare{|\DD_{b_i}| \geq n^{1-2b_{i+1}} y_{b_i}} \, .
\] 
A first estimate of the right hand term is given by Lemma 1.2 of \cite{AC04},
that we now recall.

\noindent{\bf Lemma 1.2 of \cite{AC04}}.
{\it Assume $d \geq 3$. There is a constant $\kappa_d>0$ such
that for any $\Lambda\subset \Z^d$, and any $t>0$}
\[
P\pare{ l_{\infty}(\Lambda) >t }
\leq \exp\left(-\kappa_d \frac{t}{|\Lambda|^{2/d}}\right),
\quad\text{where}\quad l_{\infty}(\Lambda)=\sum_{x \in \Lambda} l_{\infty}(x) \, .
\]
Hence, if we drop the index $i$, and set $b=b_{i+1}\approx b_i$,
for $L=n^{1-2b}y_b$, we have
\begin{eqnarray}\label{bound-often} 
\P_0 \pare{|\DD_b| \ge L}&\leq&
\sum_{\Lambda \subset ]-n;n[^d, |\Lambda|=L} \P_0 \pare{\DD_b = \Lambda,
l_n(\Lambda) \geq n^{b} L}\cr
& \leq & (2n)^{dL} \exp(-\kappa_d y_b^{1-2/d} n^{\zeta})\quad\text{with}\quad
\zeta=1-b-\frac{2}{d}(1-2b).
\end{eqnarray} 
Since $\zeta>1/2$ when $b<1/2$ and $d>4$, this estimate would suffice if
the combinatorial factor $(2n)^{dL}$
were negligible. This case corresponds to ``large'' $b$. 
For ``small'' $b$, we need to get rid of the combinatorial term.
Inspired by Le Gall's work~\cite{legall}, we propose a reduction to intersection
local times of two independent random walks. Assume indeed for a moment that 
we can compare  $\sum_{x \in \DD_b} l_n^2(x)$
with $\sum_{x \in \DD_b} l_n(x) \tilde{l}_n(x)$, where 
$(\tilde{l}_n(x))_{x \in \Z^d}$ is an independent copy of 
$(l_n(x))_{x \in \Z^d}$. Then, using Lemma 1.2 of \cite{AC04}, we obtain
\begin{eqnarray}\label{heur.eq1}
\P_0\pare{\sum_{x \in \DD_b} l^2_n(x)\geq n y_b}
& \simeq & 
\P_0\otimes\tilde\P_0 \pare{\sum_{x \in \DD_b} l_n(x) \tilde{l}_n(x) \geq n y_b}
\leq \P_0\otimes\tilde\P_0 \pare{\tilde{l}_n(\DD_b) \geq n^{1-b} y_b }\cr
& \leq 
& \E_0 \cro{\exp\pare{-\kappa_d  \frac{n^{1-b} y_b}{|\DD_b|^{2/d}}}} \cr
& \leq 
& \exp\pare{-\kappa_d n^{\zeta-2\epsilon/d} y_b}
+ \P_0 \pare{|\DD_b| \geq n^{1-2b+\epsilon}} \cr
& \leq & 
\exp\pare{-\kappa_d n^{\zeta-2\epsilon/d} y_b} 
+ \P_0 \pare{\sum_{x \in \DD_b} l^2_n(x)\geq n^{1+\epsilon}}
\end{eqnarray} 
 
Now the last term in \refeq{heur.eq1} should be negligible compared 
to the left-hand term of \refeq{heur.eq1}, so that we obtain 
\begin{equation} 
\P_0\pare{\sum_{x \in \DD_b} l^2_n(x)\geq n y_b} 
\leq \exp\pare{- c  n^{\zeta-2\epsilon/d}}  \, .
\end{equation} 

This in turns, motivates the next result, interesting on its own.
Define, for $0<b<a$, $\DD_{b,a}:=\{x:\ n^b\le l_n(x)<n^{a}\}$.
\begin{prop}\label{levelsets-prop}
Fix positive numbers $b,\gamma,\zeta$ and $\delta$ small.
The following inequality holds for large $n$
\begin{equation}\label{turb-9}
\P_0(|\DD_{b,b+\delta}|>n^\gamma y)\le \exp(-n^\zeta),\quad\text{with}\quad
\zeta<(1-\frac{2}{d})\gamma+b-\delta,
\end{equation}
provided we assume either
\begin{enumerate}
\item (i) $\gamma+2b>1$, and $y>0$, or (ii) $\gamma+2b=1$, and
$y>1+\sum_{x\in \Z^d}G_d^2(x)$.
\item $b>\frac{2}{d}\gamma$, in which case, we can take 
$\zeta=(1-\frac{2}{d})\gamma+b$.
\end{enumerate}
\end{prop}  

The paper is organized as follows. 
We gather the technical Lemmas, and the proof of Proposition
\ref{levelsets-prop} in Section~\ref{fund.sec}.
The results of Section~\ref{fund.sec} are applied to the 
problem of large deviations for SILT in Section~\ref{SILT-prop2.dem}.
We give also in Section \ref{SILT-prop2.dem}, the proof of
Proposition~\ref{lnp-prop} as well as 
large deviations estimates for $\sum_{\DD} l^p_n(x)$ for $p \neq 2$,
where $\DD$ are subsets of the range of the random walk.
 In Section~\ref{RWRS.sec}, we treat
the problem of large and moderate deviations upper bounds for the RWRS, 
and prove Proposition~\ref{propMDUB}. 
Finally, the corresponding 
lower bounds (Proposition \ref{propMDLB}) are shown in Section~\ref{RWRSLB.sec}.

%___________________________________________________________________
% 
% LEMMES FONDAMENTAUX.
%____________________________________________________________________

\section{Technical Lemmas.}
\label{fund.sec}
\subsection{Estimates for {\it low} level sets}
\begin{lemma}\label{key.lem}
Assume $d \geq 5$, and fix positive real numbers $b,\gamma,z$. Let
\[
\Db=\acc{x \in \Z^d: l_n(x) \leq n^{b} } \, .
\]
Assume that either (i) $\gamma=1$ and $z>1 + 2 \sum_{x \in \Z^d} G_d(x)^2$ or
(ii) $\gamma>1$ and $z>0$. 
Then, for $\zeta<\gamma-b-\frac{2}{d}(\gamma-2b)$ there is
a constant $c$ (depending also on $b$, $\gamma$, $z$) 
such that for $n$ large enough, 
\begin{equation}\label{key.eq}
\P_0 \pare{\sum_{x \in \Db} l_n^2(x) \geq n^{\gamma} z} \leq \exp(-c n^{\zeta}).
\end{equation}
\end{lemma}

\noindent{\bf Proof:}
We first prove the case $\gamma=1$. The case $\gamma>1$ is less delicate
and will follow the same pattern. We then indicate the necessary changes
for the case $\gamma>1$.

\noindent{\underline{ Case $\gamma=1$}}. 
The strategy is to rewrite the restricted sum of the self-intersection times
in terms of intersections of independent random walks. 
Also, we assume for simplicity that $n$ is a power of 2, $n = 2^N$; the easy
generalisation is left to the reader. First, note that, 
\begin{equation}
\label{decomp.aa1}
\sum_{x \in \Db} l_n^2(x) 
= \sum_{x \in \Db} \sum_{1\leq k,k'\leq n} \ind\{S_k=S_{k'}=x\}
\leq n+1 + 2 Z^{(0)},
\end{equation} 
\begin{equation} 
\quad\text{with}\quad
Z^{(0)}=\sum_{x \in \Db} \sum_{0 \leq k < k' \leq 2^N} \ind\{S_k=S_{k'}=x\}
\, .
\end{equation}
Now, the estimate \reff{key.eq} is equivalent to showing that
$\P_0(Z^{(0)}\ge y 2^{N}) \leq \exp(-2^{N(\zeta)})$, 
with $y>\sum_{x} G_d(x)^2$. We now bound
$Z^{(0)}\le Z^{(1)}_1 + Z^{(1)}_2 + J^{(1)}_1$, with
\[
Z^{(1)}_1= \sum_{x}\ind\{l_{2^{N-1}}(x) \leq 2^{Nb}\}
 \, \sum_{0 \leq k < k' \leq 2^{N-1}} \ind\{S_k=S_{k'}=x\},
\]
\[
Z^{(1)}_2= \sum_{x}\ind\{l_{2^N}(x)-l_{2^{N-1}}(x) \leq 2^{Nb}\} \,
\sum_{2^{N-1} \leq k < k' \leq 2^{N}} \ind\{S_k=S_{k'}=x\},
\]
\[
J^{(1)}_1= \sum_{x}\ind\{l_{2^{N-1}}(x) \leq 2^{Nb}\} \,
\sum_{0 \leq k \leq 2^{N-1} \leq  k' \leq 2^{N}} \ind\{S_k=S_{k'}=x\}.
\]
We can express $Z^{(1)}_1, Z^{(1)}_2$ and $J^{(1)}_1$ in terms of
the two independent random walks 
\[
\forall k\in \{0,\dots,2^{N-1}\}\quad
S_{k,1} = S_{2^{N-1}} - S_{2^{N-1}-k},
\quad\text{and}\quad 
S_{k,2} = S_{2^{N-1}} - S_{2^{N-1}+k}.
\]
Indeed, denoting by $\{l_{k,i}(x),k \in \N,x\in \Z^d\}$ the local times
of the random walk $(S_{k,i})_{k \in \N}$, we have
on the event $\acc{S_{2^{N-1}} =y}$, $l_{2^{N-1}}(x) =l_{2^{N-1},1}(y-x)$,
and $l_{2^N}(x)-l_{2^{N-1}}(x) = l_{2^{N-1},2}(y-x)$. 
 We obtain therefore 
$Z^{(0)}\leq Z_1^{(1)}+Z_2^{(1)}+I_1^{(1)}$ with for $i=1,2$
\[
Z^{(1)}_i  =  \sum_{y\in \Z^d} \ind\{S_{2^{N-1}} =y\}
                \sum_{x}\ind\{ l_{2^{N-1},i}(y-x) \leq 2^{Nb}\}  \,
                \sum_{0 \leq k < k' \leq 2^{N-1}}
                \ind\{S_{k,i}=S_{k',i}=y-x\} \, .
\]
Changing $x$ in $y-x$ in the second summation, we obtain for $i=1,2$
\[
Z^{(1)}_i = \sum_{x}\ind\{ l_{2^{N-1},i}(x) \leq 2^{Nb}\} \,
                \sum_{0 \leq k < k' \leq 2^{N-1}}
                \ind\{S_{k,i}=S_{k',i}=x\} \, .
\]
Finally,
\[
J^{(1)}_1 = \sum_{x}\ind\{ l_{2^{N-1},1}(x) \leq 2^{Nb}\}
                l_{2^{N-1},1}(x) l_{2^{N-1},2}(x) \, .
\]
We now denote $\{y_1,\dots,y_{N}\}$ positive reals summing up to $\bar y<y$, and
$\{b_0,\dots,b_M\}$ a regular subdivision of $[0,b]$ of mesh $\delta>0$
, such that $b_0=0$, $b_M=M\delta=b$. The precise form
of $\{y_l,b_i\}$ is given later.
From $Z^{(0)}\le Z^{(1)}_1 + Z^{(1)}_2 + J^{(1)}_1$, we deduce
\begin{equation}\label{turb-1}
\P_0(Z^{(0)}>y2^N)\le \P_0( Z^{(1)}_1 + Z^{(1)}_2>\bar y_2 2^N)+\P_0(J^{(1)}_1>y_1 2^N),
\quad\text{with}\quad \bar y_2=y_2+\dots+y_N.
\end{equation}
If we define for $i=0,\dots,M-1$
\[
\DD_{i,1}^{(1)} = \acc{x:\ 2^{Nb_i} < l_{2^{N-1},1}(x) \le 2^{Nb_{i+1}}},
\]
then, the idea is to replace $\{|\DD_{i,1}^{(1)}|$ large $\}$ by a condition on
$Z^{(1)}_1$.
Thus, we introduce
\begin{equation}\label{turb-2}
\G^{(1)}=\{\forall i=0,\dots,M-1; |\DD_{i,1}^{(1)}|\le 4 \bar y_2 2^{N(1-2b_i)}\}.
\end{equation}
The symbol $\G^{(1)}$ stands for {\it good} set at the first generation. 
Now, note that only sites visited more than once appear in
$\{\DD_{i,1}^{(1)}, i=0,\dots,M-1\}$ and contribute to $Z^{(1)}_1$.
Thus, using $k(k-1)\geq k^2/2$ for an integer $k \geq 2$, we have
\begin{eqnarray} \label{Z1.eq}
Z^{(1)}_1=\sum_{x \in \DD^{(1)}_{i,1}} \sum_{1 \leq k < k'\leq 2^{N-1}} \!\!\!
\ind\acc{S_{k,1}=S_{k',1}=x} 
& = & \frac{1}{2} \sum_{x \in \DD^{(1)}_{i,1}} 
l_{2^{N-1},1}(x) (l_{2^{N-1},1}(x)-1) \cr
& \geq & \frac{1}{4} \sum_{x \in \DD^{(1)}_{i,1}} l^2_{2^{N-1},1}(x) \ge
\frac{|\DD^{(1)}_{i,1}|}{4}2^{2b_iN}. 
\end{eqnarray}
Thus, from \reff{Z1.eq}
\begin{equation}\label{turb-3}
(\G^{(1)})^c \subset\{Z^{(1)}_1 + Z^{(1)}_2> \bar y_2 2^N\}.
\end{equation}
Thus, \reff{turb-1} becomes
\begin{equation}\label{turb-4}
\P_0(Z^{(0)}>y2^N)\le 2\P_0( Z^{(1)}_1 + Z^{(1)}_2>\bar y_2 2^N)+
\P_0(\G^{(1)},J^{(1)}_1>y_1 2^N).
\end{equation}
We proceed now by induction, and at generation $l$, we have $2^{l}$ independent
strands whose local times we denote by $l_{2^{N-l},k}$. We introduce
for $k=1,\dots,2^{l-1}$ 
\begin{equation}\label{turb-5}
\DD^{(l)}_{i,k}=\{x:\  2^{Nb_i}\le l_{2^{N-l},2k-1}(x)< 2^{Nb_{i+1}}\},
\end{equation}
and for $i=0,\dots,M-1$
\begin{equation}\label{turb-6}
J^{(l)}_{k,i} = \sum_{x}\ind\{ x\in \DD^{(l)}_{i,k}\}
                l_{2^{N-l},2k-1}(x) l_{2^{N-l},2k}(x),\quad\text{and}\quad
J^{(l)}_{k}=\sum_{i=0}^{M-1} J^{(l)}_{k,i}.
\end{equation}
The {\it good} sets at generation $l$ are as follow. We first set 
$\bar y_{l+1}=y_{l+1}+\dots+y_N$, and
\[
\forall k=1,\dots,2^{l-1},\quad\G^{(l)}_k=\{\forall i=0,\dots,M-1;
|\DD^{(l)}_{i,k}|<4 \bar y_{l+1} 2^{N(1-2b_i)}\},\quad 
\G^{(l)}=\bigcap_k \G^{(l)}_k.
\]
As in \reff{turb-3}, we obtain
\begin{equation}\label{Zl.eq}
(\G^{(l)})^c \subset 
\{Z^{(l)}_1+\dots + Z^{(l)}_{2^l}> \bar y_{l+1} 2^N\}.
\end{equation}
It is easy, after $N$ inductive steps, to obtain
\begin{equation}\label{turb-7}
\P_0(Z^{(0)}>y2^N)\le 2^N\sum_{l=1}^N \P\pare{\G^{(l)},\sum_{k=1}^{2^{l-1}}
J^{(l)}_k> y_l 2^N}, 
\end{equation}
where for each $l$, the random variables $(J^{(l)}_k, k=1,\cdots 2^{l-1})$
are i.i.d and are distributed as a variable, say $J^{(l)}$, with
\[ 
J^{(l)}=\sum_{x}1\{l_{2^{N-l}}(x) \leq 2^{Nb}\} 
l_{2^{N-l}}(x) \tilde{l}_{2^{N-l}}(x), 
\] 
where $\{\tilde{l}_{n}(x), x \in \Z^d\}$  is an independent copy of
$\{l_{n}(x), x \in \Z^d\}$. The strategy is now the following:

\begin{itemize}
\item When $l$ is {\it large}, we use the trivial bound $J^{(l)}_k
\leq 2^{2(N-l)}$,
and the classical Cramer's estimates for sums of i.i.d. In that case,
we need to center the $J^{(l)}_k$'s, i.e. to have $y_l 2^N > 2^{l-1} E[J^{(l)}]$.

\item When $l$ is {\it small}, the trivial bound $J^{(l)}_k
\leq 2^{2(N-l)}$ is to crude. To use Cramer's estimates, 
we need the existence of some exponential moments $J^{(l)}$. 
\end{itemize}

First, we specify the $\{y_l\}$. They have to be chosen in order to center
the variables $J^{(l)}$. Set $I_{\infty} = 
\sum_x l_{\infty}(x) \tilde{l}_{\infty}(x)$, where $l$ and $\tilde{l}$
are the local times of two independent walks. Note that for $d \geq 5$, 
$m_1= \E_0[I_{\infty}]  = \sum_{x \in \Z^d} G_d(x)^2 <\infty$. 
A convenient choice is the following. 
Set $l^*=N(1-\delta_0)$ ($\delta_0$ small enough), and choose 
\begin{itemize}
\item $y_l=y/(2N)$ for $l<l^*$.
\item $y_l=y/2^{N-l+2}$, for $l\geq l^*$.
\end{itemize}
It is easy to check that $\sum_{l=1}^N y_l \leq (1-\delta_0) y/2 + y/2 \leq y$.
 We obtain for large enough $N$, 
a decomposition of $\P_0(Z^{(0)}>y2^N)
\leq 2^N(R_1+R_2)$ 
\begin{equation}\label{turb-8}
R_1=\sum_{l<l^*} \P_0
\pare{\G^{(l)},\sum_{k=1}^{2^{l-1}} \bar J^{(l)}_k\geq \frac{2^N y}{4N}},
\quad\text{and}\quad R_2=\sum_{l^*\le l\le N}
\P_0\pare{\sum_{k=1}^{2^{l-1}} \bar J^{(l)}_k\geq 2^{l-1}(y-m_1)},
\end{equation}
with centered variables $\bar{J}^{(l)}_k=J_k^{(l)}-\E_0[J_k^{(l)}]$.

\noindent
\underline{About $R_2$}.
Note that for all $l$, $J^{(l)} \leq 2^{2(N-l)}$.
Using Markov inequality, for any $\lambda > 0$,
\[
\P_0 \pare{ \sum_{k=1}^{2^{l-1}} \bar{J}^{(l)}_k \geq 2^{l-1} (y-m_1)}
\leq \exp\pare{-\lambda 2^{l-1} (y-m_1) }
	\E_0 \cro{\exp\pare{ \lambda \bar{J}^{(l)}}}^{2^{l-1}} \, .
\]
We choose $\lambda = 1/2^{2(N-l)}$  and use the fact that 
$\exp(x) \leq 1 + x + 2 x^2 $ for $|x| \leq 1$, to obtain
\[
\E_0 \cro{\exp\pare{ \lambda \bar{J}^{(l)}}}
\leq 1 + 2 \lambda^2 E[(\bar{J}^{(l)})^2] \leq 1+2m_2\lambda^2 \, , 
\]
where $m_2=\E_0\cro{I_{\infty}^2}
< \infty$ by~\cite{kmss}. Thus,
\begin{equation}\label{RLDUBT2.eq}
R_2 \leq \sum_{l\geq l^*}e^{-2^{l-1}\lambda\pare{(y-m_1)- 2 m_2 \lambda}}
\leq \sum_{l \geq l^*}  e^{-2^{3l-2N-1}( y-m_1 -2m_2\lambda)}
\leq  N e^{-c 2^{N(1-3 \delta_0)}} \, .
\end{equation}  
Hence, for $\delta_0<1/(2d)$, $R_2$ is much smaller than $\exp(-2^{N\zeta})$ 
with $\zeta=1-2/d-b(1-4/d)<1-2/d<1-3\delta_0$.

\noindent
\underline{About $R_1$}.
We first obtain the existence of some exponential moments for $J^{(l)}$. 
For each $l\le l^*$, $k=1,\dots,2^{l-1}$ and any $u>0$, we use
Lemma 1.2 of \cite{AC04}, and independence between $l_{2^{N-l},2k-1}$
 and $l_{2^{N-l},2k}$,
\begin{eqnarray}\label{S2AA.3}
\P_0(J^{(l)}_k>u,\G^{(l)}_k) &\leq & \sum_{i=0}^{M-1}
\P_0\pare{\G^{(l)}_k,\sum_{\DD^{(l)}_{i,k}}
 \tilde{l}_{2^{N-l}}(x)> \frac{u}{2^{Nb_{i+1}}M}}\cr
&\leq &  \sum_{i=0}^{M-1} \E_0\cro{\ind\{|\DD^{(l)}_{i,k}|<4 y 2^{N(1-2b_i)}\}
\exp\pare{-\frac{\kappa_d u}{2^{Nb_{i+1}} |\DD^{(l)}_{i,k}|^{2/d}M}}}\cr
&\leq &\sum_{i=0}^{M-1}
 \exp\pare{-\frac{C }{2^{N\zeta_i}M}u}\leq
M\exp \pare{-\frac{C }{2^{N\max_i(\zeta_i)}M}u} \, ,
\end{eqnarray}
with $\zeta_i=b_{i+1}+\frac{2}{d}(1-2b_i)=i\delta(1-4/d)+2/d+\delta$.
Thus,
\[
\max_i(\zeta_i) = M\delta(1-4/d)+2/d+\delta \leq b (1-4/d)+2/d+\delta,
\]
and for any $\epsilon>0$, we can choose $\delta$ such that
$\max(\zeta_i) \leq b+\frac{2}{d}(1-2b)+\frac{\epsilon}{2}$.
Thus, we have a constant $C$ such that
\begin{equation}\label{AC.eq}
\P_0(J^{(l)}_k>u,\G^{(l)}_k)\leq \exp(-\xi_N u),
\quad \text{with}\quad \xi_N=\frac{C }{M} 2^{-N(b+(1-2b)2/d+\epsilon/2)}.
\end{equation}
Note that, when $u> \kappa_s/\xi_N^2$,
this estimate is better than an estimate obtained from \cite{kmss},
\begin{equation}
\label{KMSS.eq}
\P_0(J^{(l)}_k>u,\G^{(l)}_k) \leq \P_0(I_{\infty} > u) \leq \exp(-\kappa_s \sqrt{u}).
\end{equation}
However, it permits us to consider exponential moment $E[\exp(\lambda J_k)]$
for $\lambda<\xi_N$. We now go back to the standard Cramer's method. For
simplicity of notations, we drop the indices $l$ and $k$ when unambiguous.
Returning now to evaluating $R_1$, for any $0\leq \lambda<\xi_N$,
\begin{equation}
\label{S2AA.6}
\P_0 \pare{\sum_{k=1}^{2^{l-1}} \bar{J}_k^{(l)}
\geq \frac{2^{N} y}{4N};\bigcap_k \G^{(l)}_k}
\leq \exp \pare{- \lambda \frac{2^{N}y }{4N}}
\E_0 \cro{ \exp( \lambda \bar{J}^{(l)});\G^{(l)}_k}^{2^{l-1}}\, .
\end{equation}
Now, using $e^x \leq 1+x +2x^2$ for $x \leq 1$, and $\E_0\cro{\bar{J}}=0$, we have
\begin{eqnarray*}
\label{S2AA.7}
\E_0[e^{\lambda \bar{J}};\G^{(l)}_k]
&=& \E_0[e^{\lambda \bar{J}} ;\{J<1/\lambda\}]
+\E_0[e^{\lambda \bar{J}} ;\{J\geq 1/\lambda\}\cap\G^{(l)}_k]
\\[.2cm]
& \leq & \E_0[e^{\lambda \bar{J}} ;\{ J<1/\lambda\}]
+ \E_0[e^{\lambda J} ;\{J\geq 1/\lambda\}\cap\G^{(l)}_k]
\\[.2cm]
& \leq & \E_0 \cro{ \pare{1+\lambda \bar{J}+2\lambda^2(\bar{J})^2};\{J<1/\lambda\}}
+\E_0\cro{e^{\lambda J} ;\{J\geq 1/\lambda\}\cap\G^{(l)}_k}
\\[.2cm]
&\leq & 1 + \lambda \E_0\cro{|\bar{J}| ;\{J \geq 1/\lambda\}}
+ 2 \lambda^2 \E_0\cro{\bar{J}^2}
+ \E_0\cro{e^{\lambda J} ;\{J\geq 1/\lambda\}\cap\G^{(l)}_k}.
\end{eqnarray*}
Now,
\begin{eqnarray*}
\E_0\cro{|\bar{J}| ;\{J \geq 1/\lambda\}}
\leq \E_0\cro{(\bar{J})^2}^{1/2} \P_0\pare{ J \geq 1/\lambda}^{1/2}
\leq \lambda \E_0\cro{J^2} \leq \lambda \E_0(I_{\infty}^2) \, .
\end{eqnarray*}
Note that by the results of \cite{kmss},
$\E_0(I_{\infty}^2) < \infty$. Hence, for
some constant $c$,
\[
\E_0[e^{\lambda \bar{J}};\G^{(l)}_k] \leq 1 + c\lambda^2 +
\E_0\cro{e^{\lambda J};\{J\geq 1/\lambda\}\cap\G^{(l)}_k} \, .
\]
We now show that for some constant $C$,
$E\cro{e^{\lambda J};\{J\geq 1/\lambda\}\cap\G^{(l)}_k}
\leq C\lambda^2$. We decompose
this last expectation into
\[
\E_0 \cro{e^{\lambda J};\{J\geq 1/\lambda\}\cap\G^{(l)}_k}
=e^1 \P_0(\lambda J\geq 1)+I\le e^1\E_0[I_{\infty}^2] \lambda^2 +I,
\]
with
\[
I=  \int_{1/\lambda}^{\infty}\!\! \lambda e^{\lambda u}
\P_0(J^{(l)}_k\geq u;\G^{(l)}_k) du,
\quad \text{and we choose} \quad
\lambda=\frac{\xi_N}{2\log(1/\xi_N^3)}.
\]
To bound  $I$, we use estimate \refeq{AC.eq}, $\lambda<\xi_N/2$ and $N$ large enough
\begin{eqnarray}\label{S2AA.8}
I &\leq & \int_{1/\lambda}^{\infty}\!\!
 \lambda e^{\lambda u-\xi_N u}du
\leq\frac{2\lambda}{\xi_N}\int_{1/\lambda}^{\infty}\!\! (\xi_N/2)e^{-(\xi_N/2)u}du\cr
&\leq& \frac{2\lambda}{\xi_N}\exp(-\frac{\xi_N}{2\lambda})\le
\frac{\xi_N^3}{\log(1/\xi_N^3)}\le 4\xi_N \log(1/\xi_N^3) \lambda^2\le \lambda^2.
\end{eqnarray}
Thus, there is a constant $C$ such that
\[
\E_0[\exp(\lambda \bar J^{(l)}_k);\G^{(l)}_k]\leq 1+C\lambda^2\leq \exp(C\lambda^2),
\]
which together with \reff{S2AA.6}, yield
\begin{eqnarray}
\label{S2AA.9}
2^N R_1& \leq & N 2^N \exp \pare{- \frac{2^{N}y}{4N}
\frac{\xi_N}{2\log(1/\xi_N^3)} +\frac{C\xi_N^2 2^l}{4 \log^2(1/\xi_n^3)}}\cr
& \leq &
N 2^N\exp \pare{ - \frac{2^{N}y\xi_N}{16N\log(1/\xi_N^3)}} \, ,
\end{eqnarray}
where we used that  $2^{N}y> 4CN \xi_N 2^l / \log(1/\xi_N^3)$ for 
any $l\leq l^*$ and $N$ large enough, as soon as $\epsilon$ is chosen so that
$1-b-\frac{2}{d}(1-2b)-\epsilon/2 > 0$. Now, we can use
an extra $\epsilon/2$ to swallow the denominator $N\log(1/\xi_N^3)$ in
the exponential, as well as the $N2^N$ factor in front of  the exponential
in \reff{S2AA.9}.  We  obtain then for large enough $N$,
\begin{equation}
\label{S2AA.10}
\P_0(Z^{(0)}>y 2^N)
\leq \exp\pare{ - C 2^{N\zeta}},\text{ with }
\zeta=1-b-\frac{2}{d}(1-2b)- \epsilon\, .
\end{equation}
\hspace*{\fill} \rule{2mm}{2mm}

\noindent{\underline{ Case $\gamma>1$}}. 
The sequence $\{y_l\}$ and the good sets $\{\G^{(l)},l=1,\dots,N\}$ 
are different here. Since recentring the $J^{(l)}_k$ poses
no problem, we can choose $y_l=y/N$ for $l=1,\dots,N$. Next, we set for $l=1,\dots,N$
\begin{equation}\label{turb-2g}
\G^{(l)}=\{\forall i=0,\dots,M-1;
\forall k=1,\dots,2^{l-1} ; 
|\DD_{i,k}^{(l)}|\le 4\bar y_{l+1}2^{N(\gamma-2b_i)}\}.
\end{equation}
As in \reff{Zl.eq}, we obtain
\begin{equation}\label{turb-3g}
(\G^{(l)})^c
\subset\{Z^{(l)}_1+\dots + Z^{(l)}_{2^l}> \bar y_{l+1} 2^{N\gamma}\}.
\end{equation}
By induction, we obtain an inequality similar to \reff{turb-7}
with $y_l2^{N\gamma}$ replacing $y_l2^N$. 
The proof follows exactly the same pattern yielding the desired result.
\hspace*{\fill} \rule{2mm}{2mm}

%END OF KEY PROOF

\noindent{\bf Proof of Proposition~\ref{levelsets-prop}(1):}
Note that 
\[
\{|\DD_{b,b+\delta}|>n^\gamma y\}\subset \{\sum_{\underline{\DD}_{b+\delta}}l_n(x)^2
>n^{\gamma+2b}y\}.
\]
We invoke Lemma~\ref{key.lem} to conclude the proof.
\hspace*{\fill} \rule{2mm}{2mm}

\subsection{Estimates for {\it high} level sets}
We only need an improvement of Lemma 1.2 of \cite{AC04}.
\begin{lemma}
\label{fund1.lem}
Assume $d \geq 3$. There exists a constant $\kappa_d >0$ such
that for any $t > 0$, $L \geq 1$,
\begin{equation}
\label{key-amin}
\P_0 \pare{ |\acc{x: l_n(x) \geq t}| \geq L} \leq (2n)^{dL}
\exp(- \kappa_d t L^{1-2/d}) \, .
\end{equation}
\end{lemma}

\noindent{\bf Proof:}
The proof is a simple application of Lemma 1.2 of \cite{AC04}.
\begin{eqnarray}\label{key-progress}
P(|\acc{x: l_n(x) \geq t}| \geq L)
&\leq& \sum_{\Lambda \subset ]-n;n[^d; |\Lambda|=L}
P(\forall x \in \Lambda, l_n(x) \geq t) \cr
& \leq &  \sum_{\Lambda \subset ]-n;n[^d; |\Lambda|=L}
P(l_n(\Lambda) \geq Lt)\cr
& \leq&n^{dL} \exp(- \kappa_d t L^{1-2/d}).
\end{eqnarray}
\hspace*{\fill} \rule{2mm}{2mm}

\noindent{\bf Proof of Proposition~\ref{levelsets-prop}(2):}
Note that
\[
\{|\DD_{b,b+\delta}|>n^\gamma y\}\subset \{|\acc{x: l_n(x) \geq n^b}|\geq yn^\gamma\}.
\]
We use Lemma~\ref{fund1.lem} with $L=yn^\gamma$ and $t=n^b$. The
combinatorial term is negligible when $b>(2/d)\gamma$, and this gives
the correct exponent $b+(1-2/d)\gamma$.
\hspace*{\fill} \rule{2mm}{2mm}

%__________________________________________________________________
% PREUVE DE LA BORNE SUP SANS LOG.
%______________________________________________________________________

\section{Estimates for SILT.}
\label{SILT-prop2.dem}
%------------------------------------------------------
% PREUVE DE LA BORNE SUP POUR LES GRANDS TEMPS LOCAUX.
%------------------------------------------------------
We first prove Proposition \ref{LDSILT.prop}, then the lower bound
of Proposition \ref{RLDSILT.prop}, and finally estimates on
$\P_0(\sum l^p_n(x)>n^\gamma)$. 
\subsection{Proof of Proposition \ref{LDSILT.prop}}
Note that 1. of Proposition \ref{LDSILT.prop} is a direct corollary
of Lemma~\ref{key.lem}. Thus, we focus on point 2. of Proposition \ref{LDSILT.prop}.
Note first that
\begin{eqnarray*} 
\P_0\pare{\sum_{x:\, l_n(x) \geq \sqrt{n}} l_n^2(x) \geq ny}
& \leq &  
\P_0\pare{\exists x; l_n(x) \geq \sqrt{n}} 
\leq \sum_{x \in]-n;n[^d} \P_0\pare{l_n(x) \geq \sqrt{n}} 
\\
& \leq & \sum_{x \in]-n;n[^d} \P_0(H_x < \infty) \P_x(l_n(x) \geq \sqrt{n})
\\
& \leq  &
c n^d \P_0(l_n(0) \geq \sqrt{n}) \leq c n^d \exp(-\underline{c} \sqrt{n})
\, .
\end{eqnarray*} 
Thus, it is enough to prove that for any $y > 0$ and 
any $\epsilon \in ]0, 1/2 - 1/d[$, $\exists \tilde{c} > 0$ such that
\begin{equation} 
\label{grandln.eq}
\limsup_{n \rightarrow \infty} \frac{1}{\sqrt{n}} \log
\P_0\pare{\sum_{\DD_{1/2-\epsilon,1/2}} l_n^2(x) \geq n y} \leq - \tilde{c} \, ,
\end{equation} 
where for any $a,b$, with $0<b<a$, we have defined
\begin{equation}\label{defDDab}
\DD_{b,a}= \acc{x: n^{b} \leq  l_n(x) \leq n^{a}} \, .
\end{equation}
We write $\DD_{1/2-\epsilon,1/2} \subset \cup_{i=0}^{M-1} \DD_i$,  with
$b_0 \leq 1/2-\epsilon$, $b_M =1/2$. However, this 
time, $M$ will depend on $n$ (actually $M \simeq \log(\log(n)$). 
Let $(y_i, i=0 \cdots M-1)$ be  positive numbers such that
$\sum_i y_i \leq 1$.
Then, using Lemma \ref{fund1.lem},
\begin{eqnarray*} 
\P_0\pare{\sum_{\DD_{1/2-\epsilon,1/2}} l_n^2(x) \geq n y} 
 &\leq &  \sum_{i=0}^{M-1}  \P_0 \pare{\sum_{x \in \DD_i} l_n^2(x) \geq ny_i y}
\leq \sum_{i=0}^{M-1} \P_0 \pare{|\DD_i| \geq n^{1-2b_{i+1}} y_i y}\cr
&\leq& \sum_{i=0}^{M-1} n^{d n^{1-2b_{i+1}} y_i y} 
\exp \pare{-\kappa_d n^{b_i+(1-2/d)(1-2b_{i+1}) } (y_i y)^{1-2/d}} \, .
\end{eqnarray*} 
Therefore, we need to choose $(y_i, b_i, 0\leq i\leq M-1)$ such that
for some $\beta > 0$,
\begin{equation} 
\label{cond-fa2}
\left\{
\begin{array}{ll}
n^{1-2b_{i+1}} y_i \log(n) \ll n^{b_i+(1-2/d)(1-2b_{i+1}) } y_i^{1-2/d}\cr
n^{b_i+(1-2/d)(1-2b_{i+1}) } y_i^{1-2/d} \geq \beta \sqrt{n}
\end{array}
\right.
 \Leftrightarrow 
\left\{
\begin{array}{ll}
(n^{1-2b_{i+1}} y_i)^{2/d} \log(n) \ll n^{b_i}
\\
\beta n^{1/2-b_i} \leq n^{2(1-2/d)(1/2-b_{i+1}) } y_i^{1-2/d}
\end{array} \right.
\end{equation} 
For $i=M-1$, the second condition in \refeq{cond-fa2} is 
$\beta n^{1/2-b_{M-1}} \leq y_{M-1}^{1-2/d}$,
so that we have to take 
\[
1/2-b_{M-1} = 1 /\log(n),\quad\text{and}\quad
y_{M-1} = (\beta e)^{\frac{d}{d-2}}.
\]
For this choice of $b_{M-1}, y_{M-1}$, the first condition in 
\refeq{cond-fa2} is satisfied.

For the others $b_i$ ($i \leq M-2$), we take  
$b_{i+1} -1 /2 = a (b_i - 1/2)$, with $\frac{d}{2(d-2)} < a < 1$.
Hence for $i \leq M-1$, $\frac{1}{2} - b_i = (\frac{1}{a})^{M-1-i} 
\frac{1}{\log(n)}$. To have $b_0 \leq \frac{1}{2} - \epsilon < b_1$,
we take $M -1 = \lceil \frac{\log(\epsilon \log(n))}
{\log(1/a)} \rceil$. With these choices, the second condition 
in \refeq{cond-fa2} becomes 
\[
\forall i \leq M-2,\quad 
y_i \geq \beta^{\frac{d}{d-2}} \exp\pare{- 2  (1/a)^{M-i-1} (a-\frac{d}{2(d-2)})},
\]
and we take $y_i$ to satisfy the equality. Now,
the first condition in \refeq{cond-fa2} is for $i \leq M-2$,
\begin{equation} 
\beta^{\frac{2}{d-2}} \exp\pare{ \frac{d}{d-2} 
\pare{\frac{1}{a}}^{M-i-1 }} \ll \frac{\sqrt{n}}{\log(n)}
 \Leftarrow 
\beta^{\frac{2}{d-2}} \exp\pare{ \frac{d}{d-2} 
\pare{\frac{1}{a}}^{M-1 }} \ll \frac{\sqrt{n}}{\log(n)}
\, .
\end{equation} 
Recalling the value of $M$, this is satisfied as soon as 
\begin{equation} 
\label{contrainte.eps.eq}
\frac{\epsilon}{a} \pare{\frac{d}{d-2}} < \frac{1}{2} \, .
\end{equation}  
But for $\epsilon < 1/2-1/d$, one can find $a \in ]\frac{d}{2(d-2)},1[$ 
such that \refeq{contrainte.eps.eq} holds.
 
It remains now to check that we can take $\beta$ in order 
to get $\sum_{i=0}^{M-1}  y_i \leq 1$. But,
\begin{eqnarray*}
\sum_{i=0}^{M-1}  y_i 
& = & \beta^{\frac{d}{d-2}} \cro{ e^{ \frac{d}{d-2}} 
+ \sum_{i=1}^{M-1} \exp \pare{ -2 \pare{a -  \frac{d}{2(d-2)}}
\pare{\frac{1}{a}}^i } }
\\
& \leq & \beta^{\frac{d}{d-2}} \cro{ e^{ \frac{d}{d-2}} 
 + \sum_{i=1}^{\infty} \exp \pare{ -2 \pare{a - \frac{d}{2(d-2)}}
\pare{\frac{1}{a}}^i}} \, .
\end{eqnarray*} 
Since the last series is convergent, one can obviously find 
$\beta$ such that $\sum_{i=0}^{M-1} y_i \leq 1$.
\hspace*{\fill} \rule{2mm}{2mm} 

%: LOWER BOUND.
\subsection{Proof of the lower bound of Proposition \ref{RLDSILT.prop}}
For $k \in \N$, let $T^{(k)}_0$ be the $k$-th return time at $0$:
\[
T^{(0)}_0 \triangleq  0\, , \,\, 
T^{(k)}_0 \triangleq \inf \acc{n > T_0^{(k-1)}, S_n = 0} \, . \,\, 
\]
For $y > 0$,
\begin{eqnarray*}	
P\pare{\sum_{x \in \Z^d} l_n^2(x) \geq ny}
& \geq & P \pare{l_n(0) \geq \lfloor \sqrt{ny} \rfloor +1} 
 = P \pare{T_0^{(\lfloor \sqrt{ny} \rfloor)} \leq n}
\\
& \geq & P \pare{\forall k \in \{1, \cdots \lfloor \sqrt{ny} \rfloor \},
T_0^{(k)}-T_0^{(k-1)} \leq \frac{n}{\lfloor \sqrt{ny} \rfloor}}
\\
& \geq  & P \pare{T_0 \leq \frac{n}{\sqrt{ny}}}^{\sqrt{ny}} 
\\
& = & \pare{P \pare{T_0 < \infty}-
		P \pare{\frac{\sqrt{n}}{\sqrt{y}} < T_0 < \infty}  }^{\sqrt{ny}}
\end{eqnarray*} 
This proves the lower bound since 
$ \lim_{n \rightarrow \infty} P \pare{\frac{\sqrt{n}}{y} < T_0 < \infty} = 0$,
and  $P(T_0 < \infty)<1$ for $d \geq 3$.
\subsection{About $\{\sum_{\DD} l^p_n(x)>n^\gamma\}$.}
We present two results about upper bounds for $\P_0(\sum_{\DD} l^p_n(x)>n^\gamma)$,
where $\DD$ is a subset of $\{x:l_n(x)\ge 1\}$.
The first estimate concerns sites visited not too often, and is
a corollary of Lemma~\ref{key.lem}. The notation $\DD_{b,a}$ 
is defined in~\reff{defDDab}.
\begin{prop}
\label{key2.lem}
Fix positive numbers $a,b,\gamma,\zeta,p,y$, with $a>b$. 
For any $y>0$, the following inequality holds for some constant $c>0$
\begin{equation}
\label{key2.eq}
\P_0 \pare{ \sum_{x \in \DD_{b,a}} l_n^p(x) \geq n^{\gamma} y }
 \leq \exp(-c n^\zeta),
\end{equation}
provided the following conditions are satisfied.
\begin{itemize}
\item[(0)] When $p>2$, $\gamma>1+a(p-2)$. When $p\le 2$, $\gamma>1-b(2-p)$. 
\item[(i)] When $p \ge\frac{d}{d-2}$, $\zeta <\gamma-a(p-1)-\frac{2}{d}(\gamma-ap)$.
\item[(ii)] When $1<p\le\frac{d}{d-2}$, $\zeta<\gamma-b(p-1)-\frac{2}{d}(\gamma-bp)$.
\end{itemize}
The constant $c$ depends on $a,b,p,\gamma,\zeta,y$, and we take $n$ large enough.
\end{prop}

\noindent
{\bf Proof:} The strategy is to partition $\DD_{b,a}$ 
into a finite number $M$ of regions:
\begin{equation}\label{DDi}
\DD_i = \acc{x: n^{b_i} \leq l_n(x) < n^{b_{i+1}}},\quad\text{where}\quad
b_i=b+\frac{i}{M}(a-b),\quad\text{for }i=0,\dots,M,
\end{equation}
where $M$ will be chosen large enough later.
Then, using Proposition~\ref{levelsets-prop} for an arbitrarily small $\epsilon$,
one can choose $M$ large enough such that
\begin{eqnarray}
\label{ineq-fa1}
\P_0 \pare{\sum_{x \in \DD_{b,a}} l_n^p(x) \geq n^{\gamma} y}
& \leq &
\sum_{i=0}^{M-1} \P_0 \pare{\sum_{x \in \DD_i} l_n^p(x)
\geq \frac{n^{\gamma}y}{M}}\cr
& \leq &
\sum_{i=0}^{M-1} \P_0 \pare{| \DD_i| \geq \frac{n^{\gamma -pb_{i+1}}y}{M}} \cr
&\le & \sum_{i=0}^{M-1} e^{-C n^{\zeta_i}},\text{  where  }
\zeta_i=(\gamma-pb_{i+1})(1-\frac{2}{d})+b_i-\epsilon\cr
&\le & M e^{-C n^\zeta},\quad\text{where}\quad \zeta=\min_{i} \zeta_i.
\end{eqnarray}
However, Proposition~\ref{levelsets-prop} requires that $2b_i+\gamma-pb_{i+1}>1$
for each $i=0,\dots,M-1$. This is satisfied when $\gamma$ satisfies the
condition (0), and when
$(a-b)/M$ is small enough. If we set $\delta=(a-b)/M$, then we rewrite $\zeta_i$ as
\begin{equation}\label{turb-10}
\zeta_i=\gamma(1-\frac{2}{d})+b_{i+1}(1-p+\frac{2}{d}p)-\delta-\epsilon.
\end{equation}
Thus, in case (i), we have $1-p+\frac{2}{d}p\le 0$ and the minimum over 
the $\zeta_i$ is reached for the largest $b_{i+1}$. In case (ii),
the minimum over the $\zeta_i$ is reached for the smallest $b_{i+1}$. It 
then remains to choose $\delta$ small enough to obtain the desired result.
\hspace*{\fill} \rule{2mm}{2mm}

As a corollary of Lemma~\ref{fund1.lem}, we obtain the following
estimates for the regions where the local times are {\it large}.
We recall that for $p>1$, we denote by $p^*:=p/(1-p)$ the conjugate
exponent of $p$.
\begin{prop}\label{fund2.lem}
Assume $d \geq 3$, and fix positive numbers $a,b,\gamma,\zeta,p,y$, with $a>b$.
For any $y>0$, the following inequality holds for some constant $c>0$
\begin{equation}
\label{fund2.eq}
\P_0 \pare{ \sum_{x \in \DD_{b,a}} l_n^p(x)
\geq n^{\gamma} y } \leq \exp(-c n^\zeta) \, ,
\end{equation}
provided (0) $\zeta\le \frac{d}{2} b$, and either of the following two conditions.
\begin{itemize}
\item $p \geq d/(d-2)$; (i) $\zeta <\frac{\gamma }{p(2/d)+1}$;
(ii) $\zeta<\gamma-a(p-1)-\frac{2}{d}(\gamma-ap)$.
\item $1<p\le d/(d-2)$; 
(iii) $\zeta<\gamma-b(p-1)-\frac{2}{d}(\gamma-bp)$.
\end{itemize}
The constant $c$ depends on $a,b,p,\gamma,\zeta,y$, and $n$ is taken large enough.
\end{prop}

\noindent
{\bf Proof:} As in the proof of Proposition~\ref{key2.lem},
we decompose $\DD_{b,a}$ into a finite number $M$ of regions,
as in \reff{DDi}, where $M$ will be chosen later. Then, as in \reff{ineq-fa1},
\[
\P_0 \pare{\sum_{x \in \DD_{b,a}} l_n^p(x) \geq n^{\gamma} y} \leq
\sum_{i=0}^{M-1} \P_0 \pare{| \DD_i| \geq \frac{n^{\gamma -pb_{i+1}}y}{M}} \, .
\]
We now use Lemma~\ref{fund1.lem} with $t=n^{b_i}$ and
$L =n^{\gamma -pb_{i+1}}y/M$ to get
\begin{equation}
\label{region1.eq}
\P_0 \pare{\sum_{x \in \DD_{b,a}} l_n^p(x) \geq n^{\gamma}y}
\leq
\sum_{i=0}^{M-1} n^{dn^{\gamma -pb_{i+1}}y/M}
\exp\pare{- \kappa_d n^{b_i+ (1-2/d)(\gamma-pb_{i+1})} (y/M)^{1-2/d}} \, .
\end{equation}
To conclude, it is now enough to check that we can find a finite
 sequence $(b_i, 0 \leq i \leq M)$, such that $b_0= b$, $b_M > a$
and satisfying the constraints
\begin{equation}
\label{contrainte-region1}
\left\{
\begin{array}{l}
\gamma -p b_{i+1} < b_i + (1-2/d)(\gamma - p b_{i+1})
\\
\zeta \leq  b_i + (1-2/d)(\gamma -p b_{i+1})
\\
b_i < b_{i+1}
\end{array}
\right.
\Leftrightarrow
\left\{
\begin{array}{ll}
b_{i+1} > \frac{\gamma}{p} - \frac{d}{2p} b_i  &  \quad ( C_2)
\\
b_{i+1} \leq  \frac{\gamma}{p}  +  \frac{d}{p(d-2)} (b_i-\zeta)
& \quad (C_1)
\\
b_{i+1} > b_i & \quad (C_0)
\end{array}
\right.
\, .
\end{equation}
\begin{figure}[!ht]
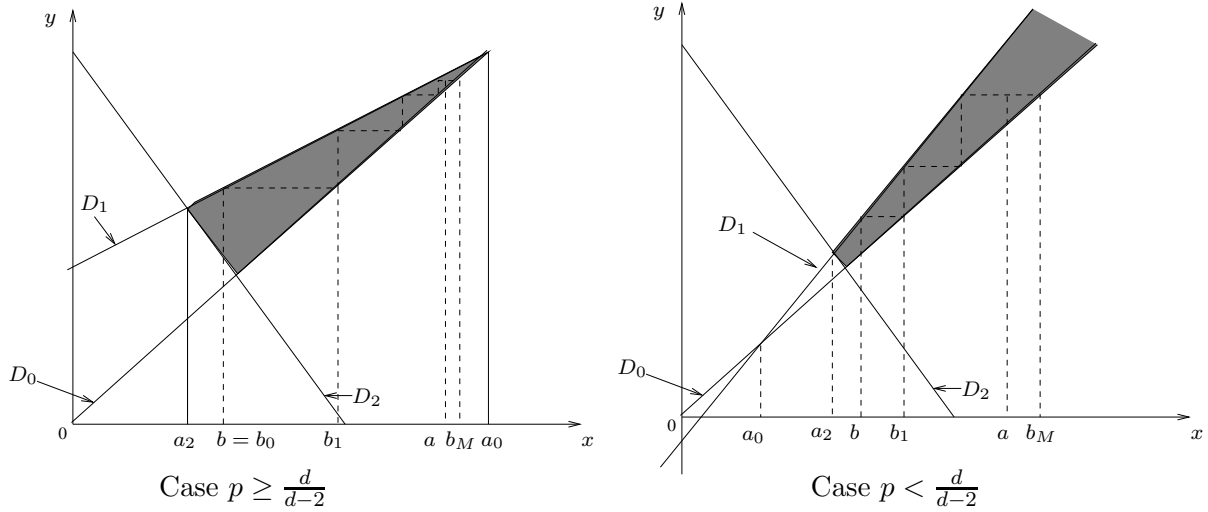

\input region1.pstex_t
\caption{Construction of $(b_i, 0 \leq i \leq M)$ for  $\DD$}
\label{region1.fig}
\end{figure}

\noindent
Let $D_0$ be the line $y=x$, $D_1$ be the line
$y=\frac{\gamma}{p}  +  \frac{d}{p(d-2)} (x-\zeta)$, and $D_2$ the
line $y=\frac{\gamma}{p} - \frac{d}{2p}x$.
Assume first that $p\not= d/(d-2)$ . Let $a_0$ (resp. $a_2$) be the abscissa of the
intersection of  $D_1$ with  $D_0$ (resp. $D_2$)
\[ a_0 = \frac{\gamma-\zeta d/(d-2)}{p-d/(d-2)} \,, \quad
a_2 = \frac{\zeta}{(d/2)} \, .
\]

\noindent
{\underline{Case $p > d/(d-2)$}:} In that case, the slope of $D_1$ is
less than 1.  Then, the region of constraints is non empty
(see figure \ref{region1.fig}) if and only if
\[
a_2 <  a_0 \Leftrightarrow \zeta < \frac{\gamma }{1+2p/d}\ (\text{i.e. condition (i)})
\, .
\]
In that case, it is always possible to construct a finite sequence
$(b_i)_{0\leq i \leq M}$ satisfying the constraints $(C_0), (C_1), (C_2)$
and $b_0=b$, $b_M \geq a$, as soon as $b \geq a_2$ i.e. (0) and
$a < a_0$, i.e. (ii). A possible choice is to take
$b_{i+1} =  \frac{\gamma}{p}  +  \frac{d}{p(d-2)} (b_i-\zeta)$, $M$
being defined by $b_{M-1} < a \leq b_M$.

\noindent
{\underline{Case $p < d/(d-2)$}:} In that case, the slope of $D_1$ is
greater than 1, and the region of constraints is never
empty.  It is always possible to construct a finite sequence
$(b_i)_{0\leq i \leq M}$ satisfying the constraints $(C_0), (C_1), (C_2)$
as soon as $b > a_0$, and $b \geq  a_2$.
 A possible choice is to take
$b_{i+1} =  \frac{\gamma}{p}  +  \frac{d}{p(d-2)} (b_i-\zeta)$, $M$
being defined by $b_{M-1} < a \leq b_M$.

When $p = d/(d-2)$, we choose $b_{i+1}-b_{i}=
\gamma/p -\zeta>0$ with the appropriate boundaries. Note that conditions
(i), (ii) and (iii) are identical and reads $\zeta <\gamma/p=\gamma(1-2/d)$.
\hspace*{\fill} \rule{2mm}{2mm}

\subsection{Proof of Proposition \ref{lnp-prop}}
\subsubsection{Point (i) of Proposition \ref{lnp-prop}}
Note first that we only need to deal with sites in $\DD_{0,\gamma/p}$.
Fix an arbitrarily small $\epsilon>0$.
We first focus on $\DD_{0,\gamma/p-\epsilon}$. We consider three cases.

\noindent{\underline{ Case $d/(d-2)<p\le 2$}}.
Proposition~\ref{key2.lem} with $b=0$ and $a=\gamma/p-\epsilon$ yields
\begin{equation}\label{lnp-eq3}
\frac{\gamma}{p}<\gamma-a(p-1)-\frac{2}{d}(\gamma-pa),
\end{equation}
since the condition (0) $\gamma>1$ holds.

\noindent{\underline{ Case $2<p$ and $\gamma> p/2$}}.
Note that Proposition~\ref{key2.lem} imposes that $a<(\gamma-1)/(p-2)$.
Note also that $\gamma>p/2$ is equivalent to $\gamma/p<(\gamma-1)/(p-2)$.
Thus, we can again take $b=0$ and $a=\gamma/p-\epsilon$ in Proposition~\ref{key2.lem},
to obtain \reff{lnp-eq3}. Condition (0) follows from $\gamma>p/2$.

\noindent{\underline{ Case $2<p$ and $\gamma\le p/2$}}.
Proposition~\ref{key2.lem} is used to deal with $\DD_{0,a}$ with $a=(\gamma-1)/(p-2)-
\epsilon$, for $\epsilon>0$ arbitrarily small. We use Proposition~\ref{fund2.lem}
to control the contribution of sites of $\DD_{a,\gamma/p-\epsilon}$. Indeed, the
three conditions we have to check reads
\begin{equation}\label{lnp-eq4}
(0)\ \frac{\gamma}{p}<\frac{d}{2}\pare{\frac{\gamma-1}{p-2}-\epsilon},\quad
(i)\ \frac{\gamma}{p}<\frac{\gamma}{2p/d+1},\quad
(ii)\ 0<\epsilon(p-1-\frac{2}{d}p).
\end{equation}
Condition (i) is equivalent to $p>d/(d-2)$ which holds here,
whereas (0) is equivalent to
$\gamma>1+2(p-2)/(4+p(d-2))$.

The proof that for some constant $C>0$
\[
\P_0(\sum_{\DD_{\gamma/p-\epsilon,\gamma/p}}
l_n^p(x)>n^\gamma)\le \exp(-C n^{\gamma/p}),
\]
is similar to the
tedious proof of Proposition~\ref{LDSILT.prop}(2), and is left to the reader.

The lower bound follows trivially from
$\{l_n(0)>n^{\gamma/p}\}\subset \{\sum l_n^p(x)>n^\gamma\}$.
\hspace*{\fill} \rule{2mm}{2mm}
\subsubsection{Point (ii) of Proposition \ref{lnp-prop}}
We assume here that $1<p\le d/(d-2)$, and $p>\gamma>1$
and we start with proving the
upper bound in \reff{lnp-eq2}. Fix an arbitrarily small $\epsilon>0$,
set $q=d/(d-2)+\epsilon$, and choose $\alpha$ such that
\[
\frac{1}{p}=\alpha+\frac{1-\alpha}{q}\Longrightarrow \alpha=1-q^*/p^*.
\]
The idea is to interpolate between $1$ and $q$. In other words,
\begin{equation}\label{lnp-eq5}
\pare{\sum_{x\in\Z^d}l^p_n(x)}^{1/p}\le 
\pare{\sum_{x\in\Z^d}l_n(x)}^{\alpha}\pare{\sum_{x\in\Z^d}l^q_n(x)}^{(1-\alpha)/q}\le
n^{\alpha}\pare{\sum_{x\in\Z^d}l^q_n(x)}^{q^*/(qp^*)}
\end{equation}
Thus, 
\[
\{\sum l^p_n(x)>n^{\gamma}\}\subset \{\sum l_n^q(x)>n^{\tilde \gamma}\},
\quad\text{with}\quad \tilde \gamma=\pare{\frac{\gamma}{p}-1+\frac{q^*}{p^*}}
\frac{qp^*}{q^*}.
\]
A simple computation yields $\gamma>1$ is equivalent to $\tilde \gamma>1$.
Thus, we can use (i) of Proposition~\ref{lnp-prop} to obtain that
\begin{equation}\label{lnp-eq6}
\zeta:=\frac{\tilde\gamma}{q}=\pare{\frac{\gamma}{p}-1+
\frac{q^*}{p^*}} \frac{p^*}{q^*}=
1-\frac{p^*(p-\gamma)}{pq^*}=1-\frac{(p-\gamma)}{(p-1)q^*}.
\end{equation}
Since this is true for any $\epsilon>0$, we have $\zeta<
1-\frac{2(p-\gamma)}{d(p-1)}$.

We prove now the lower bound in \reff{lnp-eq2}. We set $\RR_n:=\{x:l_n(x)\ge 1\}$,
and use Holder's inequality
\begin{equation}\label{lnp-eq7}
\pare{\frac{1}{|\RR_n|}\sum_{\RR_n} l_n(x)}^p\le
\frac{1}{|\RR_n|} \sum_{\RR_n} l_n(x)^p.
\end{equation}
Thus, recalling that $p>\gamma>1$, we have from \reff{lnp-eq7}
\begin{eqnarray}\label{lnp-eq8}
\P_0\pare{\sum l^p_n(x)>n^{\gamma}}&\ge&
\P_0\pare{ |\RR_n|\le r^d},\quad\text{with}\quad r^d=n^{(p-\gamma)/(p-1)}\cr
&\ge & \P_0(\sigma_r\ge n)
,\quad\text{with}\quad \sigma_r:=\inf\{k\ge 0:\ S_k\not\in ]-r/2;r/2[^d\}.
\end{eqnarray}
We use now the classical estimate $\P_0(\sigma_r\ge n)\ge \exp(-Cn/r^2)$, for
some constant $C$, and if we set $n/r^2=n^\zeta$, this yields 
$\zeta=1-\frac{2}{d}\frac{(p-\gamma)}{(p-1)}$.

\hspace*{\fill} \rule{2mm}{2mm}

%____________________________________________________________
% SECTION 5: MODERATE DEVIATION FOR RWRS. UPPER BOUNDS
%_____________________________________________________________

\section{Upper bounds for the  deviations of the RWRS}
\label{RWRS.sec}

The aim of this Section is to prove Proposition \ref{propMDUB}.
Let $\Lambda$ denote the log-Laplace transform of $\eta(0)$:
\[  \forall t \in \R\, , \Lambda(t) = \log E_{\eta} \cro{\exp(t\eta(0))} \, .
\]
Since $\eta(0)$ is centered, there exists a constant $C_0$ such that 
for $|t| \leq 1$, $\Lambda(t) \leq C_0 t^2$. By Tauberian Theorem, 
for $\eta(0)$ having the tail behavior \refeq{eq-intro.4}, 
$\Lambda(t)$ is of order $t^{\alpha^*}$ for large $t$, where
$\alpha^*$ is the conjugate exponent of $\alpha$ ($\frac{1}{\alpha}
+\frac{1}{\alpha^*} =1$). Hence, there exists a constant $C_{\infty}$ 
such that for $t \geq 1$, $\Lambda(t) \leq C_{\infty} t^{\alpha^*}$.

Our aim is to show that $P(X_n>y n^\beta)\le \exp(-C n^{\zeta})$.
In each region, we partition the range into two domains
$\Dh=\acc{x \in \Z^d; l_n(x) \geq n^{b}}$ and $\Db=
\acc{x \in \Z^d; 0<l_n(x) \leq n^{b}}$, parametrized by
a positive $b$, which will turn out to be $\beta-\zeta$.

First, for any $y_1, y_2 > 0$, such that $y_1+y_2=y$,
\begin{equation} 
\label{dec.eq}
P\pare{\sum_x \eta(x) l_n(x) \geq n^{\beta}y}
 \leq 
P\pare{\sum_{x \in \Dh} \eta(x) l_n(x) \geq n^{\beta}y_1}   
+ P\pare{\sum_{x \in \Db} \eta(x) l_n(x) \geq n^{\beta}y_2} \, .
\end{equation} 
Let $A:=\acc{\sum_{x \in \Dh} l^{\alpha^*}_n(x) \geq 
n^{\beta-b+\alpha^*b} \frac{y_1}{2C_{\infty}}}$.  
\begin{eqnarray*} 
P\pare{\sum_{x \in \Dh} \eta(x) l_n(x) \geq n^{\beta}y_1} 
& \leq &
\P_0\pare{A}
+ P\pare{ A^c \, ; 
\, \sum_{x \in \Dh} \eta(x) l_n(x) \geq n^{\beta}y_1}  
\\
& \leq &
\P_0\pare{A} + e^{-n^{\beta -b}y_1} 
	\E_0\cro{ \ind_{A^c} \exp \pare{\sum_{x \in \Dh} 
	\Lambda\pare{\frac{l_n(x)}{n^{b}}}}}
\, .
\end{eqnarray*} 
Now, on $\Dh$, $l_n(x)
\geq n^{b}$, so that using the behaviour of $\Lambda$ near infinity, 
\begin{eqnarray}\label{terme1.eq}
P\pare{\sum_{x \in \Dh} \eta(x) l_n(x) \geq n^{\beta}y_1} 
&\leq & \P_0\pare{A} + e^{-n^{\beta -b}y_1} 
\E_0\cro{ \ind_{A^c} \exp \pare{C_{\infty} 
\frac{\sum_{x \in \Dh} l_n^{\alpha^*}(x)} {n^{\alpha^* b}}}}\cr
& \leq &
\P_0\pare{\sum_{x \in \Dh} l^{\alpha^*}_n(x) \geq 
	n^{\beta-b+\alpha^*b} \frac{y_1}{2 C_{\infty}}} 
+ e^{-n^{\beta -b} y_1/2} \, .
\end{eqnarray}  
Thus, we need to prove in each region that for some constant $C>0$,
and the appropriate parameters $\beta,b$, and $\alpha$, we have
\begin{equation}\label{Reg1.eq}
\P_0\pare{\sum_{x \in \Dh} l^{\alpha^*}_n(x) \geq
        n^{\beta-b+\alpha^*b} \frac{y_1}{2 C_{\infty}}}
\le e^{-Cn^{\beta -b} }.
\end{equation}
Similarly, but using this time the behaviour of $\Lambda$
near $0$,
\begin{equation}
\label{terme2.eq}
P\pare{\sum_{x \in \Db} \eta(x) l_n(x) \geq n^{\beta}y_2} 
\leq 
\P_0\pare{\sum_{x \in \Db} l_n^2(x) \geq 
	n^{\beta+b} \frac{y_2}{2 C_0}} 
+ \exp(-n^{\beta -b} y_2/2) \, .
\end{equation}  
In this case, we need to prove that for some constant $C>0$
\begin{equation}\label{Reg2.eq}
\P_0\pare{\sum_{x \in \Db} l_n^2(x) >
        n^{\beta+b} \frac{y_2}{2 C_0}} \le \exp(-C n^{\beta -b})
\end{equation}
\vspace{.5cm}
\noindent{\bf Region I}.  
We choose $b=1-\beta$ in order to have $\zeta_I=\beta-b=2\beta-1$.
We first prove \refeq{Reg2.eq}. Since $\beta+b=1$, Lemma~\ref{key.lem} requires that 
$y_2 >y_0:=2 (1 +2 \sum_x G^2_d(x))C_0$. The condition of
Lemma~\ref{key.lem} on $\zeta_I$ reads
\begin{equation}\label{ineq-I}
\beta-b< \beta+b-b-\frac{2}{d}(\beta+b-2b)\Longrightarrow 
\beta<b(1+d/2)\Longrightarrow \beta<\frac{1+d/2}{2+d/2}.
\end{equation}
Secondly, \refeq{Reg1.eq} relies on Proposition~\ref{fund2.lem} with
$p=\alpha^*$, $\gamma=\beta+b(\alpha^*-1)$ and $\zeta=\beta-b$.
Condition (0) is equivalent to $\beta\le(1+d/2)b$, already fulfilled
in \reff{ineq-I}.

When $\alpha>d/2$, the condition (iii) of Proposition~\ref{fund2.lem}
is equivalent to $\beta<(1+d/2)b$, 
which we have already taken into account in~\reff{ineq-I}.

When $\alpha\le d/2$, condition (i) of Proposition~\ref{fund2.lem}
is also equivalent to $\beta<(1+d/2)b$.
Now, note that there is no point in considering sites visited more than
$n^{\zeta_I}$. Indeed, for some constant $C>0$
\[
\P_0 \pare{\sum_{{\bar\DD}_{\zeta_I}} l^{\alpha^*}_n(x)
\geq n^{\beta-b+\alpha^*b} y}
\leq \P_0\pare{\exists x \in ]-n;n[^d; l_n(x) \geq n^{\zeta_I}}
\leq e^{-C n^{\zeta_I}} \, .
\]
Thus, condition (ii) with $a=\zeta_I$ is equivalent to $\beta<(\alpha+1)b$,
which implies $\beta<\frac{\alpha+1}{\alpha+2}$. 

\vspace{.5cm}
\noindent 
{\bf Region II}. We choose $b=\frac{\beta}{\alpha+1}$, to get 
$\zeta_{I\!\!I}=\beta-b$. Here $\alpha<d/2$.
We start with proving \reff{Reg2.eq}. Lemma~\ref{key.lem}
imposes $\beta+b\ge 1$ and $\beta<(1+d/2)b$. The latter inequality
holds true when $\alpha<d/2$, whereas the former requires 
$\beta\ge \frac{\alpha+1}{\alpha+2}$. Note that in case $\beta+b=1$, we
need that $y_2>y_0$.

In order to prove \refeq{Reg1.eq}, we use Proposition~\ref{fund2.lem} in case
$p>(d/2)^*$ and need
to check its conditions (0),(i) and (ii). Condition (0) and (i) 
are equivalent to $\alpha\le d/2$. Finally,
Condition (ii) has to be checked with $\zeta_{I\!\!I}=\beta - b=\alpha b$ and
$\gamma=\beta - b + \alpha^* b = (\alpha + \alpha^*) b = \alpha \alpha^* b$. 
If we choose $a=\zeta_{I\!\!I}$, a simple computation yields $\zeta_{I\!\!I}=\gamma-
a(\alpha^*-1)-2/d(\gamma-a\alpha^*)$. Thus, Proposition~\ref{fund2.lem}
allows to conclude that for any $\epsilon>0$,
 \[
\P_0 \pare{\sum_{\DD_{b,\alpha b-\epsilon}} l^{\alpha^*}_n(x)
\geq n^{\alpha \alpha^*b} y} \leq \exp(-C n^{\zeta_{I\!\!I}}) \, .
\]
Hence, it remains to prove that for $y>0$, $\epsilon>0$, and 
$n$ sufficiently large,
\[
\P_0 \pare{\sum_{\DD_{\alpha b -\epsilon,\alpha b}} l^{\alpha^*}_n(x)
\geq n^{\alpha \alpha^*b} y} \leq \exp(-C n^{\zeta_{I\!\!I}}) \, .
\]
We are in the situation of point 2. of Proposition
\ref{LDSILT.prop}. The proof is the same, and is left to the reader.

We now prove \reff{major-RegII}.
We need to show that $\underline{\DD}_{b-\delta}$ and ${\bar\DD}_{b+\delta}$
bring a negligible contribution.
If we define for each $\delta>0$, $B_{\delta}:=\acc{\sum_{x \in 
\underline{\DD}_{b-\delta}} l^{2}_n(x) \geq n^{\beta-b} \frac{y_2}{2C_0}}$, then
as in \reff{terme2.eq}, we obtain 
\[
P(\sum_{\underline{\DD}_{b-\delta}}\eta(x) l_n(x) > y_2n^\beta)
\le P(B_{\delta})+e^{-Cn^{\beta-b+\delta}},
\]
and we need to show that $P(B_{\delta})\le \exp(-Cn^{\beta-b+\delta'})$ for
some $\delta'>0$. By Lemma~\ref{key.lem}, we need $\delta$ small enough
so that $\beta+b-\delta>1$. We also need to check that 
\[
\zeta_{I\!\!I}-\delta'<\gamma-\delta-(b-\delta)-\frac{2}{d}\pare{
\gamma-\delta-2(b-\delta)}\Longrightarrow
\delta+\frac{d}{2}\delta'<\frac{d/2-\alpha}{\alpha+1} \beta.
\]
Now, for the large level sets, let $A_{\delta}:=
\acc{\sum_{x \in {\bar\DD}_{b+\delta}} l^{\alpha^*}_n(x) \geq
n^{\beta-b+\alpha^*b+(\alpha^*-1)\delta} \frac{y_1}{2C_{\infty}}}$.
As in \reff{terme1.eq}, we obtain
\begin{equation}\label{level-deltaII}
P(\sum_{{\bar\DD}_{b+\delta}} \eta(x) l_n(x)> y_1 n^\beta)\le 
\P_0(A_{\delta})+e^{-Cn^{\beta-b+\delta}},
\end{equation}
and we need to show that $P(A_{\delta})\le \exp(-Cn^{\beta-b+\delta'})$ for
some $\delta'>0$. We invoke again Proposition~\ref{fund2.lem}
in the case $\alpha<d/2$.
Here, condition (0) imposes $b(d/2-\alpha)\ge\delta'-(d/2)\delta$. Condition
(i) imposes
\[
b\alpha^*(1-\frac{\alpha}{d/2})>\delta'(1+\frac{\alpha^*}{d/2})-(\alpha^*-1)\delta.
\]
Finally, condition (ii) yields $(\alpha^*-1)\delta> \alpha^* \delta'$.
Thus, conditions (0),(i) and (ii) are clearly satisfied for $\delta$ and $\delta'$
small enough.

\vspace{.5cm}
\noindent{\bf Region III}. We choose $b=2\beta/(d+2)$, and obtain 
$\zeta_{I\!\!I\!\!I}=\beta-b$. Here, $\alpha\ge d/2$, and
we need to prove a result weaker than \reff{Reg2.eq}. For any $\epsilon>0$
\begin{equation}\label{Reg3.eq}
\P_0\pare{\sum_{x \in \Db} l_n^2(x) >
        n^{\beta+b} \frac{y_2}{2 C_0}} \le \exp(- n^{\zeta_{I\!\!I\!\!I}-\epsilon})
\end{equation}
This is a direct application of Lemma~\ref{key.lem}, as soon as we check
that $\gamma=\beta+b\ge 1$ (and $y_2>y_0$ in the case
of equality). This last condition means that $\beta\ge
\frac{d/2+1}{d/2+2}$ which defines precisely Region III.

Now, we prove \reff{Reg1.eq} invoking Proposition~\ref{fund2.lem} with
$p=\alpha^*<(d/2)^*$. Condition (0) is equivalent to $b(1+d/2)\ge \beta$, and
we have here equality. Condition (iii) 
holds for any $\zeta=\beta-b-\epsilon$ by an straightforward computation.

We prove now \reff{major-RegIII}. We define $A_{\delta}$ as in Region II, and
\reff{level-deltaII} follows similarly. We invoke again Proposition~\ref{fund2.lem}
in the case $\alpha\ge d/2$ with $\gamma=\beta-(\alpha^*-1)(b+\delta)$, and
$\zeta=\beta-b+\delta'$. Condition (0) requires $(d/2) \delta\ge \delta'$,
whereas condition (iii) requires $(2/d)\delta> \delta'$. So that a choice
$\delta'=\delta/d$ yields \reff{major-RegIII}, for any $\delta>0$.
\hspace*{\fill} \rule{2mm}{2mm} 

%---------------------------------------------------------------
% SECTION 6: MODERATE DEVIATIONS FOR RWRS. LOWER BOUNDS.
%----------------------------------------------------------------

\section{Lower Bounds for RWRS.}
\label{RWRSLB.sec}

This Section is devoted to the proof of Proposition \ref{propMDLB}. The
symmetry assumption simplifies the proof, thanks to the following Lemma.

\begin{lemma} (Lemma  2.1 of \cite{AC04})
When $\{\eta(x),x\in \Z^d\}$ are independent and have bell-shaped densities,
then for any $\Lambda$ finite subset of $\Z^d$, and any $y>0$
\begin{equation}
\label{ineq-handy}
P\left(\sum_{x\in\Lambda} \alpha_x \eta(x)>y\right)\le
P\left(\sum_{x\in\Lambda} \beta_x \eta(x)>y\right),
\quad\text{if}\quad 0\leq \alpha _x\le \beta_x\text{ for all }x\in \Lambda.
\end{equation} 
\end{lemma} 

\noindent
{\bf Region I}. Let $\RR_n:= \acc{x: l_n(x) \geq 1}$.
Under the symmetry assumption, $\forall c >0$, 
\[
P\pare{\sum_x \eta(x) l_n(x) \geq n^{\beta}y} 
\geq P\pare{\sum_{x \in \RR_n} \eta(x) \geq n^{\beta}y}
\geq \P_0(|\RR_n| \geq cn) P_{\eta}\pare{\sum_{j=1}^{cn} \eta_j \geq 
n^{\beta} y} \, .
\] 
Now, it is well known, that for $d \geq 3$, there is
$c >0$ such that $\lim_{n \rightarrow \infty} \P_0(|\RR_n| \geq cn) =1$.
For the other terms, if $1/2 < \beta < 1$, we are in a regime
of moderate deviations for a sum of i.i.d., and there is $C>0$ such that 
\[
\liminf_{n \rightarrow \infty} \frac{1}{n^{2\beta-1}} 
\log P_{\eta}\pare{\sum_{j=1}^{cn} \eta_j \geq n^{\beta} y} \geq -C.
\]
\vspace{.5cm}
\noindent{\bf Region II}.
Under the symmetry assumption,
\[
P\pare{\sum_x \eta(x) l_n(x) \geq n^{\beta}y}  
\geq P\pare{\eta(0) l_n(0) \geq n^{\beta}y}
\geq P_{\eta} \pare{\eta(0) \geq n^{\frac{\beta}{\alpha+1}} y}
\P_0 \pare{l_n(0) \geq n^{\frac{\beta \alpha }{\alpha +1}}} \, .
\]
Now, for $\frac{\beta \alpha }{\alpha +1} \leq 1$, the second probability 
is of order $\exp(-C n^{\frac{\beta \alpha }{\alpha +1}})$, which
is also the order of the first one. 
This leads to the lower bound in region II.

\vspace{.5cm}\noindent{\bf Region III}.
We keep the notations of the heuristic discussion of Region III:
$T=n^\beta$, and $r^d=T^{d/(d+2)}$. Recall that $\RR_n$ is the range
of the walk, and let $\sigma_r:=\inf\{k\ge 0:\ S_k\not\in ]-r/2;r/2[^d\}$.
Under the symmetry assumption, for any $\epsilon>0$
\begin{eqnarray}\label{DV-lb1}
P\pare{\sum_x \eta(x) l_n(x)\geq yn^{\beta}}
&\geq&P\pare{\sum_x \eta(x) l_T(x)\geq yT}\cr
&\ge &P\pare{\{\forall x\in \RR_T,\ \eta(x)>y\}\cap \{\epsilon r^d<|\RR_T|<r^d\}}\cr
&\ge& P_{\eta}\pare{ \eta(0)>y}^{\epsilon r^d}\pare{\P_0(|\RR_T|<r^d)-
\P_0(|\RR_T|<\epsilon r^d)}\cr
&\ge & P_{\eta}\pare{ \eta(0)>y}^{\epsilon r^d}\pare{\P_0(\sigma_r>T)-
\P_0(|\RR_T|<\epsilon T^{d/(d+2)})}.
\end{eqnarray}
It is now well known that
for some constant $C>0$, $\P_0(\sigma_r>T)\ge \exp(-CT/r^2)$. On the other hand,
from Donsker and Varadhan ~\cite{DV}, there is a constant $c_{DV}$ such that
for $\lambda>0$
\begin{eqnarray}\label{DV-lb2}
\P_0(|\RR_T|<\epsilon T^{d/(d+2)})
&\le&\exp(\lambda \epsilon T^{d/(d+2)}) \E_0\cro{e^{-\lambda |\RR_T|}}\cr
&\le& \exp(-(c_{DV}\lambda^{2/(d+2)}-\lambda \epsilon )T^{d/(d+2)})\cr
&\le & \exp(-\frac{c_{DV}}{2}(\frac{c_{DV}}{2\epsilon})^{2/d}T^{d/(d+2)}),
\end{eqnarray}
where we have chosen $c_{DV}\lambda^{2/(d+2)}=2\lambda \epsilon$. Thus, we can
choose $\epsilon $ small enough so that $2\P_0(|\RR_T|<\epsilon r^d)\le
\P_0(|\RR_T|<r^d)$, and conclude the lower bound.

\hspace*{\fill} \rule{2mm}{2mm}

%____________________________________________________________________
% BIBLIOGRAPHIE
%___________________________________________________________________

\end{document}